\documentclass[a4paper,11pt,leqno,twoside]{article}

\usepackage{amsmath, amsthm, amssymb,bm}
\usepackage{graphics}

\makeatletter
\newif\ify@autoscale \y@autoscaletrue \def\Yautoscale#1{\ifnum #1=0
  \y@autoscalefalse\else\y@autoscaletrue\fi}
\newdimen\y@b@xdim
\newdimen\y@boxdim \y@boxdim=13pt
\def\Yboxdim#1{\y@autoscalefalse\y@boxdim=#1}
\newdimen\y@linethick    \y@linethick=.3pt
\def\Ylinethick#1{\y@linethick=#1}
\newskip\y@interspace \y@interspace=0ex plus 0.3ex
\def\Yinterspace#1{\y@interspace=#1}
\newif\ify@vcenter   \y@vcenterfalse
\def\Yvcentermath#1{\ifnum #1=0 \y@vcenterfalse\else\y@vcentertrue\fi}
\newif\ify@stdtext   \y@stdtextfalse
\def\Ystdtext#1{\ifnum #1=0 \y@stdtextfalse\else\y@stdtexttrue\fi}
\newif\ify@enable@skew   \y@enable@skewfalse
\y@vcentertrue
\def\y@vr{\vrule height0.8\y@b@xdim width\y@linethick depth 0.2\y@b@xdim}
\def\y@emptybox{\y@vr\hbox to \y@b@xdim{\hfil}}
\ify@enable@skew
 \def\y@abcbox#1{\if :#1\else
   \y@vr\hbox to \y@b@xdim{\hfil#1\hfil}\fi}
 \def\y@mathabcbox#1{\if :#1\else
   \y@vr\hbox to \y@b@xdim{\hfil$#1$\hfil}\fi}
\else
 \def\y@abcbox#1{\y@vr\hbox to \y@b@xdim{\hfil#1\hfil}}
 \def\y@mathabcbox#1{\y@vr\hbox to \y@b@xdim{\hfil$#1$\hfil}}
\fi
\def\y@setdim{%
  \ify@autoscale%
   \ifvoid1\else\typeout{Package youngtab: box1 not free! Expect an
     error!}\fi%
   \setbox1=\hbox{A}\y@b@xdim=1.6\ht1 \setbox1=\hbox{}\box1%
  \else\y@b@xdim=\y@boxdim \advance\y@b@xdim by -2\y@linethick
  \fi}
\newcount\y@counter
\newif\ify@islastarg
\def\y@lastargtest#1,#2 {\if\space #2 \y@islastargtrue
  \else\y@islastargfalse\fi}
\def\y@emptyboxes#1{\y@counter=#1\loop\ifnum\y@counter>0
  \advance\y@counter by -1 \y@emptybox\repeat}
\def\y@nelineemptyboxes#1{%
  \vbox{%
    \hrule height\y@linethick%
    \hbox{\y@emptyboxes{#1}\y@vr}
    \hrule height\y@linethick}\vspace{-\y@linethick}}
\def\yng(#1){%
  \y@setdim%
  \hspace{\y@interspace}%
  \ifmmode\ify@vcenter\vcenter\fi\fi{%
  \y@lastargtest#1,
  \vbox{\offinterlineskip
    \ify@islastarg
     \y@nelineemptyboxes{#1}
    \else
     \y@ungempty(#1)
    \fi}}\hspace{\y@interspace}}
\def\y@ungempty(#1,#2){%
  \y@nelineemptyboxes{#1}
  \y@lastargtest#2,
  \ify@islastarg
   \y@nelineemptyboxes{#2}
  \else
   \y@ungempty(#2)
  \fi}
\def\y@nelettertest#1#2. {\if\space #2 \y@islastargtrue
  \else\y@islastargfalse\fi}
\def\y@abcboxes#1#2.{%
  \ify@stdtext\y@abcbox#1\else\y@mathabcbox#1\fi%
  \y@nelettertest #2.
  \ify@islastarg\unskip%
   \ify@stdtext\y@abcbox{#2}\else\y@mathabcbox{#2}\fi%
  \else\y@abcboxes#2.\fi}
\ify@enable@skew
 \newdimen\y@full@b@xdim
 \newcount\y@m@veright@cnt
 \def\y@get@m@veright@cnt#1#2.{%
   \if :#1 \advance\y@m@veright@cnt by 1\y@get@m@veright@cnt#2.\fi}
 \let\y@setdim@=\y@setdim
 \def\y@setdim{%
   \y@setdim@ \y@full@b@xdim=\y@b@xdim
   \advance\y@full@b@xdim by 1\y@linethick}
 \def\y@m@veright@ifskew#1{
   \y@m@veright@cnt=0 \y@get@m@veright@cnt#1.
   \moveright \y@m@veright@cnt\y@full@b@xdim}
\else
 \def\y@m@veright@ifskew#1{}
\fi
\def\y@nelineabcboxes#1{%
  \y@nelettertest #1.
  \ify@islastarg
   \y@m@veright@ifskew{#1}
    \vbox{
      \hrule height\y@linethick%
      \hbox{\ify@stdtext\y@abcbox#1\else\y@mathabcbox#1\fi\y@vr}
      \hrule height\y@linethick}\vspace{-\y@linethick}
  \else
   \y@m@veright@ifskew{#1}
    \vbox{
      \hrule height\y@linethick%
      \hbox{\y@abcboxes #1.\y@vr}%
      \hrule height\y@linethick}\vspace{-\y@linethick}
  \fi}
\def\young(#1){%
  \y@setdim%
  \hspace{\y@interspace}%
  \y@lastargtest#1,
  \ifmmode\ify@vcenter\vcenter\fi\fi{%
  \vbox{\offinterlineskip
    \ify@islastarg\y@nelineabcboxes{#1}%
    \else\y@ungabc(#1)%
    \fi}}\hspace{\y@interspace}}
\def\y@ungabc(#1,#2){%
  \y@nelineabcboxes{#1}%
  \y@lastargtest#2,
  \ify@islastarg\y@nelineabcboxes{#2}%
  \else\y@ungabc(#2)%
  \fi}
\makeatother

\numberwithin{equation}{section}

\newcommand{\E}{\mathbb{E}}
\renewcommand{\S}{\mathbf{S}}

\newcommand{\Z}{\mathcal{Z}}

\newcommand{\U}{\mathbf{U}}

\newcommand{\R}{\mathbb{R}}
\newcommand{\C}{\mathbb{C}}

\newcommand{\cc}{\mathfrak{c}}
\newcommand{\RC}{\operatorname{RC}}

\newcommand{\wt}{\operatorname{wt}}
\newcommand{\Cat}{\operatorname{Cat}}

\newcommand{\supp}{\operatorname{supp}}

\def\({ \left( }
\def\){ \right)}

\theoremstyle{plain}
\newtheorem{theorem}{Theorem}[section]
\newtheorem{proposition}[theorem]{Proposition}
\newtheorem{lemma}[theorem]{Lemma}
\newtheorem{corollary}[theorem]{Corollary}
\theoremstyle{definition}
\newtheorem{definition}{Definition}[section]
\newtheorem{example}{Example}[section]

\theoremstyle{conjecture}

\theoremstyle{problem}

\date{\empty}

\title{\bfseries Jucys-Murphy Elements and Unitary Matrix Integrals}
\author{\textsc{Sho MATSUMOTO} \\
{\it Graduate School of Mathematics, Nagoya University.} \\
{\it  \small Furocho, Chikusa-ku, Nagoya, 
464-8602, JAPAN.} \\
\texttt{sho-matsumoto@math.nagoya-u.ac.jp} \\
and \\
\textsc{Jonathan NOVAK} \\
{\it Department of Mathematics, Massachusetts Institute of Technology} \\
{\it  \small Cambridge MA 02139, USA.} \\
\texttt{jnovak@math.mit.edu} \\
}

\begin{document}

\maketitle

\begin{abstract}
	In this paper, we study the relationship between polynomial integrals on the unitary group and the conjugacy class
	expansion of symmetric functions in Jucys-Murphy elements.  Our main result is an explicit formula for the top
	coefficients in the class expansion of monomial symmetric functions in Jucys-Murphy elements, from which we recover
	the first order asymptotics of polynomial integrals over $\U(N)$ as $N \rightarrow \infty$.
	Our results on class expansion include an analogue of Macdonald's result for the top connection 
	coefficients of the class algebra, a generalization of Stanley and Olshanski's result on the
	polynomiality of content statistics on Plancherel-random partitions,
	and an exact formula for the multiplicity of the class of full cycles in the expansion of a complete symmetric function in Jucys-Murphy
	elements.  The latter leads to a new combinatorial interpretation of the Carlitz-Riordan central factorial numbers.
\end{abstract}

\section{Introduction}

\subsection{Polynomial integrals on the unitary group}
	Let $\U(N)$ denote the group of $N \times N$ unitary matrices, equipped with the normalized Haar measure, and let
	$\mathcal{A} \subset L^2(\U(N),\text{Borel}, \text{Haar})$ denote the algebra of polynomial functions on $\U(N)$.  Thus 
	$\mathcal{A}$ is the set of functions
	$f: \U(N) \rightarrow \C$ for which there exists a polynomial $p_f$ in $N^2$
	variables with $f(U)=p_f(u_{11},\dots,u_{NN})$ for all $U = (u_{ij}) \in \U(N)$.  How can one compute inner products in $\mathcal{A}$?  
		
	In addition to being a natural question concerning one of the classical groups, this problem is of importance in the theory
	of random matrices and various related branches of mathematics and 
	mathematical physics, see e.g. \cite{BB, Berkolaiko, Collins, CS, DGa, DH, KMS, Rains, Sam, Xu}.  
	By linearity of the integral, the problem reduces to the computation of monomial inner products 
	
		\begin{equation*}
			\langle u_{i(1)j(1)} \dots u_{i(m)j(m)}, u_{i'(1)j'(1)} \dots u_{i'(n)j'(n)} \rangle_N
		\end{equation*}
		
	\noindent
	where $i,j:[m] \rightarrow [N]$ and $i',j':[n] \rightarrow [N]$ are given functions.
	Using the invariance of the Haar
	measure, one may show that this inner product vanishes unless $m=n$.  Thus $\mathcal{A}$ admits the orthogonal decomposition
	
		\begin{equation*}
			\mathcal{A} = \bigoplus_{n=0}^{\infty} \mathcal{A}^{(n)},
		\end{equation*}

	\noindent
	where $\mathcal{A}^{(n)}$ is the space of homogeneous polynomial functions of degree $n$.
	This orthogonality is equivalent to the statement that the the space $V^{\otimes m} \otimes V^{*\otimes n}$, where
	$V$ is the defining representation of $\U(N)$ and $V^*$ is the dual representation, admits a non-trivial space of invariants
	if and only if $m=n$.  When $m=n \leq N$, the space of $\U(N)$-invariants in the representation $V^{\otimes m} \otimes V^{*\otimes n}$
	has a basis indexed by the elements of the symmetric group $\S(n)$ \cite{Howe}.  This fact can be used \cite{Collins, CS} to derive the formula
	
		\begin{equation}	
			\label{eqn:CUEWick}
			\langle u_{i(1)j(1)} \dots u_{i(n)j(n)}, u_{i'(1)j'(1)} \dots u_{i'(n)j'(n)} \rangle_N = \sum_{\pi \in \S(n)} 
			\bigg{(} \sum_{\sigma\rho^{-1} = \pi} \delta_{i,i'\sigma} \delta_{j,j'\rho} \bigg{)} \langle u_{11} \dots u_{nn}, u_{1\pi(1)} \dots u_{n\pi(n)} \rangle_N,
		\end{equation}	
		
	\noindent
	valid for any $n \leq N$, which further reduces the problem of evaluating inner products in $\mathcal{A}$ to the evaluation of certain canonical
	inner products labelled by permutations. 
	
	The formula \eqref{eqn:CUEWick} is reminiscent of the Wick formula \cite{Janson}, which reduces the computation of mixed moments of Gaussian random 
	variables to the computation of covariances modulo the application of a combinatorial rule.  However, 
	the Wick-type rule \eqref{eqn:CUEWick} does not reduce the problem to the computation of covariances, but only to the computation of
	the basic inner products 
	$\langle u_{11} \dots u_{nn}, u_{1\pi(1)} \dots u_{n\pi(n)} \rangle_N$ indexed by permutations.  These basic inner products are
	are certain rational functions of $N$ of which depend on the cycle structure of $\pi$ in a rather complicated way.
	Nevertheless, it is known that the behaviour of these averages simplifies dramatically in the limit $N \rightarrow \infty$.  Given a permutation
	$\pi \in \S(n)$, recall that its \emph{reduced cycle type} is the partition $\mu$ whose parts $\mu_1 \geq \mu_2 \geq \dots$ are the lengths of the 
	cycles of $\pi$, each reduced by one in order to remove fixed points \cite{Macdonald}.  For example, the reduced cycle type
	of the permutation $\pi=(1\ 2\ 3\ 4)(5\ 6)(7)(8) \in \S(8)$ is $\mu=(3,1)$.  Let $C_\mu(n)$ denote the conjugacy class of permutations in $\S(n)$
	of reduced cycle type $\mu$.  Then, for any $\pi \in C_\mu(n)$,
	
		\begin{equation}
			\label{eqn:Asymptotics}
			(-1)^{|\mu|}N^{n+|\mu|} \langle u_{11} \dots u_{nn}, u_{1\pi(1)} \dots u_{n\pi(n)} \rangle_N = \prod_{i=1}^{\ell(\mu)} \Cat_{\mu_i}
			+O\bigg{(} \frac{1}{N^2} \bigg{)}
		\end{equation}
		
	\noindent
	as $N \rightarrow \infty$, where
	
		\begin{equation*}
			\Cat_r=\frac{1}{r+1}{2r \choose r}
		\end{equation*}
		
	\noindent
	are the Catalan numbers.  This estimate has been known to physicists for some time, see e.g. \cite{BB} for an application to quantum transport
	in mesoscopic systems.  
	
	The first proof of \eqref{eqn:Asymptotics} is due to
	Collins \cite{Collins}, and is quite involved.  The appearance of the Catalan numbers,
	which are ubiquitous in enumerative combinatorics \cite[Exercise 6.19]{Stanley:EC2}, motivates the search for 
	combinatorial structure underlying the averages $\langle u_{11} \dots u_{nn},u_{1\pi(1)} \dots u_{n\pi(n)} \rangle_N$ which might
	lead to a conceptually simple proof of \eqref{eqn:Asymptotics}.  
	Along these lines, it was observed in \cite{Novak:BC} that the evaluation of the basic inner products is intimately related to a certain 
	``inverse problem'' in the character theory of the symmetric groups.  This connection is explored in the present article.
	Working on the symmetric group side, we obtain a number of results concerning this ``inverse problem'' which give a new perspective on the structure 
	of polynomial integrals on the unitary group (leading in particular to a new proof of \eqref{eqn:Asymptotics}) 
	and are moreover of significant interest in their own right.
	
	\subsection{The class expansion problem}
	Consider the group algebra $\C[\S(n)]$ of the symmetric group, and let $\Z(n)$ be its center.  Associate to the conjugacy class
	$C_\mu(n)$ its indicator function, i.e. the group algebra element
		
		\begin{equation*}
			\cc_\mu(n) := \sum_{\pi \in C_\mu(n)} \pi.
		\end{equation*}
		
	\noindent
	The set $\{\cc_\mu(n) : \wt(\mu) \leq n\}$, where $\wt(\mu):=|\mu|+\ell(\mu)$, is 
	a natural basis of $\Z(n)$, which for this reason is often referred to as the \emph{class algebra}.  
	The weight statistic $\wt(\mu)$ is simply the size of the support of any $\pi \in C_\mu(n)$, i.e. the number of
	of points on which $\pi$ acts non-trivially, while $|\mu|$ is the minimal number of transpositions required to 
	generate $\pi$.
	
	Recall that the equivalence  classes of complex, irreducible representations $\mathsf{V}^\lambda$
	of $\C[\S(n)]$ are indexed by partitions $\lambda$ of size $n$, see e.g. \cite{OV}.  By Schur's Lemma, the central element
	$\cc_\mu(n)$ acts as a scalar operator in any irreducible representation $\mathsf{V}^\lambda$,
	
		\begin{equation*}
			\cc_\mu(n) \cdot \mathsf{v} = \omega_\mu(\lambda) \mathsf{v}
		\end{equation*}	
		
	\noindent
	for all $\mathsf{v} \in \mathsf{V}^\lambda$.
	The scalar $\omega_\mu(\lambda)$ is called the \emph{central character} of $\cc_\mu(n)$ acting in $\mathsf{V}^\lambda$.
			
	A fundamental problem in the representation theory of the symmetric groups is to describe the central character $\omega_\mu(\lambda)$,
	viewed as a function of $\lambda$ with $\mu$ fixed, in terms
	of natural statistics on $\lambda$.  For example, a well-known formula
	going back to Frobenius asserts that the central character of the conjugacy class $\cc_{(1)}(n)$ of transpositions acting in 
	$\mathsf{V}^\lambda$ is simply
	
		\begin{equation}
			\label{eqn:Frobenius}
			\omega_{(1)}(\lambda) = \sum_{\Box \in \lambda} c(\Box),
		\end{equation}
		
	\noindent
	the sum of the contents of the Young diagram $\lambda$.  Analogous formulas are known for more complicated partitions $\mu$,
	but their structure is very complex, see \cite{Lassalle:explicit} and references therein.
		
	The ``inverse problem'' referred to above is the following: construct a group algebra element in $\Z(n)$ whose central character in 
	$\mathsf{V}^\lambda$ admits a simple description.  One such construction is the Okounkov-Pandharipande theory of completed
	cycles \cite{OP}, which exhibits remarkable connections with Gromow-Witten theory.  We now describe an older and closely
	related construction, originally due to Jucys \cite{Jucys}, which turns out to be intimately related to the computation of
	the basic inner products $\langle u_{11} \dots u_{nn}, u_{1\pi(1)} \dots u_{n\pi(n)} \rangle_N$ in the algebra $\mathcal{A}$.
		
	The \emph{Jucys-Murphy elements}
	$J_2,\dots,J_n \in \C[\S(n)]$ are the transposition sums
	
		\begin{equation*}
			\begin{split}
				J_2 &= (1\ 2) \\
				J_3 &= (1\ 3) + (2\ 3)\\
				J_4 &= (1\ 4) + (2\ 4) + (3\ 4) \\
					&\vdots \\
				J_n &= (1\ n) + (2\ n) + \dots + (n-1\ n).
			\end{split}
		\end{equation*}
		
	\noindent
	These elements commute, and in fact $\C[J_2,\dots,J_n]$ is a maximal commutative subalgebra of $\C[\S(n)]$ known as the 
	\emph{Gelfand-Zetlin algebra}, see \cite{OV}.  Although the Jucys-Murphy elements are clearly non-central, Jucys observed that
	symmetric functions of them \emph{are} central.
	More precisely, let $\Lambda$ denote the $\C$-algebra of symmetric functions \cite{Macdonald} and 
	consider the multiset (or ``alphabet'') $\Xi_n=\{\{J_2,\dots,J_n,0,0,\dots\}\}$.
	Jucys proved that
	
		\begin{equation}
			\label{eqn:JMCentral}
			f(\Xi_n) = f(J_2,\dots,J_n,0,0,\dots) \in \Z(n)
		\end{equation}
		
	\noindent
	for any symmetric function $f$.
	We thus have a homomorphism of commutative $\C$-algebras $\Lambda \rightarrow \Z(n)$ defined by $f \mapsto f(\Xi_n)$, which will
	be called the \emph{Jucys-Murphy specialization}.  
	
	Jucys proved the validity of \eqref{eqn:JMCentral} by showing directly that
	the image of the elementary symmetric function $e_k$ under the JM specialization is
	
		\begin{equation}
			\label{eqn:JMElementary}
			e_k(\Xi_n) = \sum_{|\mu|=k} \cc_\mu(n),
		\end{equation}
		
	\noindent
	which may be visualized as the sum of all permutations on the $k^{\text{th}}$ level of the Cayley graph of $\S(n)$,
	generated by the set of all transpositions.
	For a proof of this fact, see \cite{DG, Jucys} or Proposition \ref{prop:Jucys} below.
	Since $\Lambda = \C[e_1,e_2,\dots]$ is the algebra of polynomials in the elementary symmetric functions (the ``fundamental theorem'' of
	symmetric function theory), \eqref{eqn:JMElementary} implies \eqref{eqn:JMCentral}.  Moreover, it is a classical result of Farahat and
	Higman \cite{FH} that the class sums \eqref{eqn:JMElementary} generate the algebra $\Z(n)$, so that the JM specialization is in 
	fact a surjection onto the class algebra.  
	
	Jucys' second remarkable discovery concerning the eponymous elements $\Xi_n$ 
	is that the central character of $f(\Xi_n)$ acting in $\mathsf{V}^\lambda$ is 
	given by the simple substitution rule $f(\Xi_n) \mapsto f(A_\lambda)$, where $A_\lambda = \{\{c(\Box) : \Box \in \lambda\}\}$ is the 
	content alphabet of $\lambda$.  That is, we have the action
	
		\begin{equation}
			\label{eqn:JucysCharacter}
			f(\Xi_n) \cdot \mathsf{v} = f(A_\lambda) \mathsf{v}
		\end{equation}
		
	\noindent
	for all $f \in \Lambda, \mathsf{v} \in \mathsf{V}^\lambda$.  A proof of this property may be found in \cite{Okounkov:JM}.
	Thus the JM specialization furnishes group algebra elements whose central character in a given representation takes an especially simple form, namely
	a symmetric polynomial function of the contents of the partition indexing the representation.
	As an example, one can conclude that the central character of the Farahat-Higman generator \eqref{eqn:JMElementary}
	acting in $\mathsf{V}^{\lambda}$ is simply 
	
		\begin{equation*}
			e_k(A_\lambda) = \sum_{1 \leq i_1 <  i_2 < \dots < i_k \leq n} c(\Box_{i_1})c(\Box_{i_2}) \dots c(\Box_{i_k}),
		\end{equation*}
		
	\noindent
	where $\Box_1,\Box_2,\dots,\Box_n$ is an arbitrary enumeration of the cells of the Young diagram $\lambda$.
	When $k=1$, this reduces to Frobenius' formula \eqref{eqn:Frobenius}.
	More generally, given a symmetric function $f \in \Lambda$, if one can explicitly resolve the central element $f(\Xi_n)$ into a linear
	sum of conjugacy classes $\cc_\mu(n)$, then one will have explicitly constructed an element of the group algebra whose central
	character in the irreducible representation $\mathsf{V}^\lambda$ is $f(A_\lambda)$.  This is  
	the \emph{class expansion problem}: given $f \in \Lambda$, explicitly determine its coordinates
	
		\begin{equation*}
			f(\Xi_n) = \sum_{\mu} G_\mu(f,n) \cc_\mu(n)
		\end{equation*}
		
	\noindent
	relative to the conjugacy class basis.  
	
	The class expansion problem has been completely solved for two polynomial bases of $\Lambda$, namely the elementary 
	symmetric functions $\{e_k\}$ and the power-sum symmetric functions $\{p_k\}$ \cite{LT}.  Ideally, one would like to obtain the solution 
	for a \emph{linear} basis of $\Lambda$.  In this article, we focus on the class expansion of the monomial symmetric functions $\{m_\lambda\}$,
	which constitute the most natural linear basis of $\Lambda$.  This choice of basis
	is further motivated by the surprising fact that its solution would enable one to integrate arbitrary polynomial functions on $\U(N)$;
	as observed in \cite{Novak:BC}, we have
	
		\begin{equation}
			\label{eqn:fundamental}
			\langle u_{11} \dots u_{nn},u_{1\pi(1)} \dots u_{n\pi(n)}\rangle_N = \frac{1}{N^n} \sum_{k=0}^{\infty} (-1)^k\frac{G_\mu(h_k,n)}{N^k}
		\end{equation}
		
	\noindent
	for any $\pi \in C_\mu(n)$ and $N \geq n$, where $h_k$ is the complete homogeneous symmetric function of degree $k$.  Heuristically, this formula
	can be thought of as follows.  Since the columns of a unitary matrix are unit vectors, one expects that the entries of a Haar-distributed
	random unitary matrix should be of order $1/N$, and indeed it is true that $\langle u_{ij},u_{ij} \rangle_N = 1/N$.  Furthermore,
	for any fixed $n$, any $n \times n$ submatrix of an $N \times N$ random Haar-distributed unitary matrix converges to a matrix of 
	independent random variables as $N \rightarrow \infty$ (see \cite{Novak:IMRN,PR} for a precise statement),
	so the inner product $\langle u_{11} \dots u_{nn}, u_{11} \dots u_{nn} \rangle$
	should be of order $N^{-n}$.  The (convergent) power series in \eqref{eqn:fundamental} gives the corrections to this estimate for
	$N \geq n$ finite and $\pi \in \S(n)$ different from the identity.

	\subsection{Main results and organization}
	In Section \ref{sec:MonomialTop} we present our main result:
	an explicit formula for the class coefficients $G_\mu(m_\lambda,n)$ in the ``top'' case $|\mu|=|\lambda|$.
	This formula involves certain refinements of the Catalan numbers originally considered by Haiman \cite{Haiman}.  Via the fundamental relation
	\eqref{eqn:fundamental}, the main formula yields a conceptually simple derivation of the first-order asymptotics \eqref{eqn:Asymptotics} as a corollary.	
	
	In Section \ref{sec:Connection}, we obtain general properties of the class coefficients $G_\mu(f,n)$ analogous to the properties of 
	the connection coefficients governing multiplication of conjugacy classes in the algebra $\Z(n)$.  We prove that the class coefficients
	$G_\mu(f,n)$ are polynomial functions of $n$, and as an application of this result we generalize a theorem of Stanley \cite{Stanley:content} and 
	Olshanski \cite{Olshanski} on the polynomiality of content statistics on Plancherel-random partitions.
	We furthermore use Lagrange inversion to give an analogue of Macdonald's 
	construction \cite{GJ:top, Macdonald} for top connection coefficients
	in the setting of top class coefficients.
	
	In Section \ref{sec:CentralFactorial}, we use the character theory of the symmetric groups to obain 
	an exact formula for the multiplicity of the class of full cycles
	in the conjugacy class expansion of a complete homogeneous symmetric function in Jucys-Murphy elements.
	Quite surprisingly, the formula obtained involves the central factorial numbers introduced by Carlitz and Riordan \cite{CR}
	in their study of the divided difference operator.  One consequence of this result is an elegant new combinatorial
	interpretation of the central factorial numbers as counting ``primitive'' factorizations of a full cycle into transpositions.
	Another consequence is an exact formula for cyclic inner products in $\mathcal{A}$ previously stated by Collins \cite{Collins}.
	
	Finally, in the Appendix we give explicit examples of the class expansion of monomial symmetric functions as well as 
	tables of top class coefficients.
	
	\subsection{Acknowledgements}
	Both authors would like to thank Beno\^it Collins for helpful and inspiring conversations.  Parts of this paper were written while 
	the second author was visiting 	Ecole Normale 
	Sup\'erieure de Lyon, and J. N. would like to thank Alice Guionnet for the opportunity to visit.
	
	An earlier version of this paper led to further recent research into the class expansion problem and its relationship
	with group integrals, see \cite{Feray, Lassalle:JM,Matsumoto,Zinn}.
	We would like to acknowledge helpful
	conversations and correspondence with Valentin F\'eray and Michel Lassalle regarding their works \cite{Feray, Lassalle:JM}.

\section{Top class coefficients of $m_\lambda(\Xi_n)$}
	\label{sec:MonomialTop}
	
	The monomial symmetric functions are perhaps the most natural class of symmetric functions; they are constructed by taking
	an arbitrary monomial as a seed, and growing a power series from this monomial by symmetrizing it.
	More precisely,
	let $x_1,x_2,x_3,\dots$ be formal variables, and let $\lambda=(\lambda_1 \geq \lambda_2 \geq \lambda_3 \geq \dots)$ be 
	a partition.  By definition, the monomial symmetric function of type $\lambda$ in the variables $x_i$ is the formal power series
	
		\begin{equation*}
			m_\lambda(x_1,x_2,x_3,\dots) = \sum x_1^{a_1}x_2^{a_2}x_3^{a_3} \dots,
		\end{equation*}
		
	\noindent
	where the summation is over all distinct permutations $(a_1,a_2,a_3,\dots)$ of the vector $(\lambda_1,\lambda_2,\lambda_3,\dots)$.
	Since the latter vector has only $\ell(\lambda)$ non-zero entries, this definition yields a well-defined formal power series each of 
	whose terms is a monomial of degree $|\lambda|=\lambda_1+\lambda_2+\lambda_3+\dots$.  
	If we set $x_{n+1}=x_{n+2}= \dots =0$, the monomial symmetric
	function specializes to a symmetric polynomial $m_\lambda(x_1,x_2,\dots,x_n)$ which is homogenous of degree $n$ provided 
	$\ell(\lambda) \leq n$, but is identically zero if $\ell(\lambda)>n$.
	Many other basic families of symmetric functions may be obtained from the 
	monomial symmetric functions, for example the elementary symmetric functions $e_k=m_{(1^k)}$, the Newton power-sums 
	$p_k=m_{(k)}$, and the complete symmetric functions $h_k=\sum_{\lambda \vdash k} m_\lambda$.  General references on 
	symmetric function theory are \cite{Macdonald} and \cite[Chapter 7]{Stanley:EC2}.
	
	In this section, we study the image of the monomial symmetric functions in the JM specialization.  Let $L_\mu^\lambda(n):=G_\mu(m_\lambda,n)$,
	so that
	
		\begin{equation*}
			m_\lambda(\Xi_n) = \sum_\mu L_\mu^\lambda(n)\cc_\mu(n).
		\end{equation*}
		
	\noindent
	Since the alphabet $\Xi_n$ contains $n-1$ non-zero elements, we may restrict our study to partitions $\lambda$ with $\ell(\lambda) \leq n-1$.
	From the definition of $m_\lambda$ and the JM elements, one finds that the class coefficient $L_\mu^\lambda(n)$ has the following 
	combinatorial interpretation: it is equal to the number of factorizations of a fixed representative $\pi$ of the conjugacy class $C_\mu(n)$
	into $|\lambda|$ transpositions of the form
	
		\begin{equation}
			\label{eqn:PrimitiveFactorization}
			\pi = \underbrace{(*\ 2) \dots (*\ 2)}_{a_1} \underbrace{(*\ 3) \dots (*\ 3)}_{a_2} \dots \underbrace{(*\ n) \dots (*\ n)}_{a_{n-1}},
		\end{equation}
		
	\noindent
	where $(a_1,a_2,\dots,a_{n-1})$ is a permutation of $(\lambda_1,\lambda_2,\dots,\lambda_{n-1})$.  
	Since the minimal number of transpositions required to generate a permutation of reduced cycle type $\mu$ is $|\mu|$,
	it is clear from this combinatorial interpretation that $L_\mu^\lambda(n)$ vanishes unless $|\mu| \leq |\lambda|$.  Furthermore,
	in order for $L_\mu^\lambda(n)$ to be non-zero, we must have $|\mu|+2g=|\lambda|$ for some integer $g \geq 0$; this
	is simply the fact that any permutation is either even or odd (Corollary \ref{cor:MinimalParity} below).  In this section 
	we explicitly evaluate the class coefficient $L_\mu^\lambda(n)$ in the ``top'' case $|\mu|=|\lambda|$ by counting
	all factorizations \eqref{eqn:PrimitiveFactorization} with the minimal number of factors.
	
	\subsection{Jucys' theorem}
	Before treating the class expansion problem for arbitrary monomial symmetric functions, let us revisit the case $m_{(1^k)}=e_k$ 
	considered by Jucys.
		
	\begin{proposition}[\cite{Jucys}]
			\label{prop:Jucys}
			For any $k \geq 0,$
				\begin{equation*} 
					e_k(\Xi_n) =\sum_{|\mu| = k} \cc_\mu(n).
				\end{equation*}
		\end{proposition}
		
		\begin{proof}
			Consider first the case $k\geq n.$  Then it is clear
			that $e_k(\Xi_n)=0,$ since each term in the sum 
			
				\begin{equation*}
					e_k(\Xi_n) = \sum_{2 \leq t_1 \leq \dots \leq t_k} J_{t_1} \dots J_{t_k}
				\end{equation*}
				
			\noindent
			defining $e_k(\Xi_n)$
			is a product
			of $k$ distinct factors and the alphabet $\Xi_n$ only contains $n-1$ non-zero
			elements.  Moreover, if $|\mu| \geq n$ then $\wt(\mu) \geq n+1,$ and hence
			the conjugacy class $C_{\mu}(n)$ is empty and $\cc_{\mu}(n)=0,$ whence
			the sum on the right hand side of the desired equality is $0.$  Thus the claim holds
			for $k \geq n.$
						
			Now suppose $0 \leq k \leq n-1.$  Then
				$$\sum_{|\mu|=k} \cc_{\mu}(n) = \sum_{\substack{\sigma \in \S(n)\\
				\#(\sigma) = n-k}} \sigma,$$
			\noindent
			where $\#(\sigma)$ denotes the number of cycles of $\sigma \in \S(n)$.
			To prove that the claim holds in the
			range $0 \leq k \leq n-1$ we proceed by induction on $n.$
			
			If $n=2$, then $e_1(\Xi_2)= J_2= (1,2)$ and the claim is trivial.

			Let $n >2$ and suppose the claim
			holds true for $e_{k}(\Xi_{n-1})$ with any $k$. 
			We define the projection $P_n$ from $\S(n)$ to $\S(n-1)$ by
				$$
				P_n(\sigma)(i) = \begin{cases}
				\sigma(i) & \text{if $\sigma(i)\not=n$} \\
				\sigma(n) & \text{if $\sigma(i)=n$}, 
				\end{cases}
				$$
				for $\sigma \in \S(n)$ and  $1 \le i \le n-1$.
				In other words,  $P_n(\sigma)$ is defined to be the permutation whose
				cycle decomposition is obtained by 
				erasing the letter $n$ in the cycle decomposition of $\sigma$.
				For each $\tau \in \S(n-1)$, we have $P_n^{-1}(\tau)=
				\{\tau(s,n) \ | \ 1 \le s \le n-1\} \cup \{\tau \cdot (n)\}$.
				Here $\tau \cdot (n)$ is an image under the natural injection 
				$\S(n-1) \hookrightarrow \S(n)$.
				Observe $\#(\tau(s,n))=\#(\tau)$ and $\#(\tau \cdot (n))=\#(\tau )+1$. 
				Thus, the sum over permutations with exactly $n-k$ cycles equals
				$$
					\sum_{\tau \in S(n-1)} \sum_{\begin{subarray}{c} \sigma \in P_n^{-1}(\tau) \\ 					
					\#(\sigma)=n-k 
					\end{subarray}} \sigma  
					= \sum_{\begin{subarray}{c} \tau \in S(n-1) \\ \#(\tau)=n-1-k 					
					\end{subarray}} \tau \cdot (n)
					+ \sum_{\begin{subarray}{c} \tau \in S(n-1) \\ \#(\tau)=n-k 					
					\end{subarray}} \sum_{s=1}^{n-1}
					\tau (s\ n).
				$$
				By the induction hypothesis, the first sum on the right hand side
				equals $e_{k}(J_1,\dots,J_{n-1})$
				Since $e_k(x_1,\dots,x_n)=e_{k}(x_1,\dots,x_{n-1})+e_{k-1}(x_1,\dots,x_{n-1}) 				
				x_n$,
				we obtain the desired equality for $n$.
				\end{proof}
	
	\subsection{Explicit formula for top class coefficients}
	Let $\Cat_r= \frac{1}{r+1} \binom{2r}{r}$ be the $r$th Catalan number:
	
		\begin{center}
			\begin{tabular}{|c||c|c|c|c|c|c|c|c|c|} \hline 
			$r$ & $0$ & $1$ & $2$ & $3$ & $4$ & $5$ & $6$ & $7$ & $8$  \\ \hline
			$\Cat_r$  & $1$ & $1$ & $2$ & $5$ & $14$ & $42$ & $132$ & $429$ & $1430$ \\ \hline 
			\end{tabular} \ .
		\end{center}
		
	\noindent	
	It is well known that 
	Catalan numbers satisfy the recurrence 
	
		\begin{equation} \label{eq:CatalanRecurrence}
			\Cat_r= \sum_{q=0}^{r-1} \Cat_q \Cat_{r-1-q}.
		\end{equation} 
		
	\noindent
	The Catalan numbers admit a vast array of combinatorial interpretations, see \cite[Exercise 6.19]{Stanley:EC2}.
	
	We will use the following interpretation of the Catalan numbers.
	For a positive integer $k$,
	let $\E(k)$ be the set of all weakly increasing sequences $(i_1,\dots,i_k)$
	of $k$ positive integers satisfying $i_p \ge p$ for $1 \le p \le k-1$ and $i_k=k$.
	For example,
	$$\E(3)= \{(123),(133),(223),(233),(333)\}.$$
	Then, as proved below,
	the cardinality of $\E(k)$ is $\Cat_k$.
	
	Let $(i_1,\dots,i_k)$ be a weakly increasing sequence of $k$ positive integers.
	We say that $(i_1,\dots,i_k)$ is of type $\lambda \vdash k$
	if $\lambda=(\lambda_1,\lambda_2,\dots)$ is a permutation of $(b_1,b_2,\dots)$,
	where, for each $p \ge 1$, $b_p$ is the multiplicity of $p$ in $(i_1,\dots,i_k)$.

	\begin{example}
		The sequences $(1233)$, $(1334)$, and $(1134)$ are of type $(2,1,1)$,
		while the sequences $(444477799)$, $(555669999)$ are of type $(4,3,2)$.
	\end{example}

	\begin{definition}
	Given a partition $\lambda \vdash k$, the \emph{refined Catalan number}
	$\RC(\lambda)$ counts
	sequences $(i_1,\dots,i_k)$ in $\E(k)$
	of type $\lambda$.
	If $\lambda$ is the empty partition, set $\RC(\lambda)=1$.
	\end{definition}
	
	\begin{example}
	The four sequences 
	$(1444), (2444), (3444), (3334)$ in $\E(4)$ are all of type $(3,1)$,
	and indeed $\RC(3,1)=4$.
	We have $\RC(k)=\RC(1^k)=1$.
	\end{example}

	\begin{proposition} \label{prop:SumRC}
	The sum of $\RC(\lambda)$ over $\lambda \vdash k$ equals $\Cat_k$:
	$$
	\sum_{\lambda \vdash k} \RC(\lambda)=\Cat_k.
	$$
	\end{proposition}

	\begin{proof}
	This is a direct consequence of the fact that $|\E(k)|=\Cat_k$.
	\end{proof}

	\begin{example} \label{ExampleListRC}
	We give some examples of $\RC(\lambda)$ for small $|\lambda|$.
	``SUM'' stands for the sum $\sum_{\lambda \vdash k} \RC(\lambda)
	=\Cat_k$.
	\begin{center}
	\begin{tabular}{|c||c|} \hline 
	$\lambda$ & $1$  \\ \hline
	$\RC(\lambda)$ & 1  \\ \hline
	\end{tabular} \quad
	\begin{tabular}{|c||c|c||c|} \hline 
	$\lambda$ &  $2$ & $1^2$ & SUM \\ \hline
	$\RC(\lambda)$ & 1 & 1 & 2 \\ \hline
	\end{tabular} \quad
	\begin{tabular}{|c||c|c|c||c|} \hline 
	$\lambda$ &  $3$ & $21$ & $1^3$ & \text{SUM} \\ \hline
	$\RC(\lambda)$ & 1 & 3 & 1 & 5 \\ \hline
	\end{tabular}
	\end{center}
	\begin{center}
	\begin{tabular}{|c||c|c|c|c|c||c|} \hline 
	$\lambda$ & $4$ & $31$ & $2^2$ & $2 1^2$ & $1^4$ & SUM \\ \hline
	$\RC(\lambda)$ & 1 & 4 & 2 & 6 &  1 & 14 \\ \hline
	\end{tabular} 
	\end{center}
	\begin{center}
	\begin{tabular}{|c||c|c|c|c|c|c|c||c|} \hline 
	$\lambda$ & $5$ & $41$ & $32$ & $31^2$ & $2^21$ & $21^3$ & $1^5$ & SUM \\ \hline
	$\RC(\lambda)$ & 1 & 5 & 5 & 10 &  10 & 10 & 1 & 42 \\ \hline
	\end{tabular}
	\end{center}
	\end{example}

	An explicit formula for $\RC(\lambda)$ is known and given as follows.
	See \cite{Stanley:park} and also \cite{Haiman}.

	\begin{proposition}[\cite{Stanley:park}]
	For any parition $\lambda$,
	$$
	\RC(\lambda) = \frac{|\lambda|!}
	{(|\lambda|-\ell(\lambda)+1)! \, \prod_{i \ge 1} m_{i}(\lambda)!}
	= \frac{1}{|\lambda|+1} m_{\lambda}(1^{|\lambda|+1}).
	$$
	Here $m_i(\lambda)$ is the multiplicity of $i$ in $\lambda=(\lambda_1,\lambda_2,\dots)$.
	\end{proposition}

	Note that $\RC(a^m)= \frac{1}{(a-1)m+1} \binom{am}{m}$ 
	is often called a \emph{higher Catalan number} or sometimes a \emph{Fuss-Catalan
	number}.  In particular, $\RC(2^m)=\Cat_m$.
	
	\begin{definition}
		\label{definition:refinement}
		Given two partitions $\lambda, \mu \vdash k$,
		let $\R(\lambda,\mu)$ denote the set of sequences of partitions $(\lambda^{(1)}, \lambda^{(2)},\dots, \lambda^{(\ell(\mu))})$ such
		that $\lambda^{(i)} \vdash \mu_i$ for $1 \leq i \leq \ell(\mu)$ and $\lambda^{(1)} \cup \lambda^{(2)} \cup \cdots \cup \lambda^{(\ell(\mu))}=\lambda.$
		
		\noindent
		Here $\lambda^{(1)} \cup \lambda^{(2)} \cup \cdots \cup \lambda^{(\ell(\mu))}$ 
		is the partition obtained by rearranging the juxtaposed sequence of
		parts of the partitions $\lambda^{(1)},
		\lambda^{(2)},\dots, \lambda^{(\ell(\mu))}$ in weakly decreasing order.  
		If $\R(\lambda,\mu) \neq \emptyset,$
		we say that $\lambda$ is a \emph{refinement} of $\mu.$
	\end{definition}
	
	Informally, $\lambda$ is a refinement of $\mu$ if it is obtained by ``splitting'' parts of $\mu$ into smaller pieces.
	
	The following assertions are immediate:
	\begin{itemize}
	\item $\R(\lambda, (k))$ consists of one element $(\lambda)$;
	
		\smallskip
		
	\item $\R(\lambda,(1^k))$ consists of one element $(\underbrace{(1),(1), \dots,(1)}_k)$ 
	if $\lambda=(1^k)$, or is empty otherwise; 
	
		\smallskip
		
	\item $\R((k), \mu)$ consists of one element $((k))$ if $\mu=(k)$,
	and is empty otherwise;
	
		\smallskip
		
	\item $\R((1^k), \mu)$ consists of one element  $((1^{\mu_1}), (1^{\mu_2}),\dots, 
	(1^{\ell(\mu)}))$;
	
		\smallskip
		
	\item Suppose $\ell(\lambda)=\ell(\mu)$. Then 
	$\R(\lambda,\mu)$ consists of one element $((\lambda_1),(\lambda_2),\dots,(\lambda_{\ell(\lambda)}))$
	if $\lambda=\mu$, and is empty otherwise.
	
		\smallskip
		
	\item $\R(\lambda,\mu) = \emptyset$ unless $\lambda \le \mu$.
	Here $\le$ stands for the dominance partial ordering:
	$\lambda \le \mu \Leftrightarrow \lambda_1+\cdots+\lambda_i \le \mu_1+\cdots+\mu_i \
	\text{for all $i \ge 1$}$.
	(This is may be found in \cite[I (6.10)]{Macdonald}.)
	\end{itemize}

	\begin{example} \label{exa:mfR}
	The set $\R((3,2,2,1), (5,3))$ consists of two elements given by
	$((3,2),(2,1))$ and $((2,2,1),(3))$.
	\end{example}
	
	We are now ready to state our formula for top class coefficients.
	
	\begin{theorem} \label{thm:monomial}
		Let  $\mu,\lambda$ be partitions, $|\mu|=|\lambda|.$
		Then the top class coefficient $L_{\mu}^{\lambda}(n)=G_{\mu}(m_{\lambda},n)$ is given by
		
			\begin{equation*}
				L^\lambda_\mu(n)=\sum_{(\lambda^{(1)},\lambda^{(2)},\dots,\lambda^{(\ell(\mu))}) \in 
				\R(\lambda,\mu)}
				\RC(\lambda^{(1)})\RC(\lambda^{(2)}) \cdots \RC(\lambda^{(\ell(\mu))}).
			\end{equation*}
			
	\noindent		
	In particular, $L^\lambda_\mu=L^\lambda_\mu(n)$ is independent of $n$, and
	$L^\lambda_\mu$ is zero unless $\lambda$ is a refinement of $\mu$.
	\end{theorem}

	Observe that for $\lambda,\mu \vdash k$,
	$$
	L^\lambda_{(k)}=\RC(\lambda), \qquad L^\lambda_{(1^k)}=\delta_{\lambda, (1^k)},
	\qquad L^{(k)}_\mu= \delta_{\mu, (k)}, \qquad L^{(1^k)}_\mu=1, \qquad
	L^\lambda_\lambda=1.
	$$
	The equality $L^{(1^k)}_\mu=1$ is compatible with Proposition \ref{prop:Jucys}.

	Since $L_\mu^\lambda = 0$ unless $\lambda$ is a refinement of $\mu$,
	the matrix $(L^\lambda_\mu)_{\lambda, \mu \vdash k}$ is
	strictly lower unitriangular in the sense of \cite[I-6]{Macdonald}. 

	The proof of Theorem \ref{thm:monomial} in the next subsections.
	The numbers $L^\lambda_\mu$ for $\lambda, \mu \vdash k$ for $k \le 7$ are tabulated 
	in the Appendix.  Before moving on to the proof, let us state two consequences of this theorem.
	
	Define $F_\mu^k(n):=G_\mu(h_k,n)$, where $h_k = \sum_{\lambda \vdash k} m_\lambda$ is the complete symmetric
	function of degree $k$.  Thus
	
		\begin{equation*}
			h_k(\Xi_n) = \sum_\mu F_\mu^k(n) \cc_\mu(n).
		\end{equation*}
		
	\begin{corollary}
		\label{cor:complete}
		We have
		\begin{equation*}
			F_\mu^{|\mu|}=\prod_{i \ge 1} \Cat_{\mu_i}.
		\end{equation*}
	\end{corollary}

	\begin{proof}
		For $k=|\mu|,$ we have
		$$
		F_\mu^k = \sum_{\lambda \vdash k} L^\lambda_\mu
		=\sum_{\lambda \vdash k} \sum_{(\lambda^{(1)},\lambda^{(2)},\dots) \in 
		\R(\lambda,\mu)} \RC(\lambda^{(1)})\RC(\lambda^{(2)}) \cdots \RC(\lambda^{(\ell(\mu))})
		$$
		by Theorem \ref{thm:monomial}. By the definition of $\R(\lambda,\mu)$,
		we see that
		$$
		\bigsqcup_{\lambda \vdash k} \R(\lambda,\mu) =
		\{(\lambda^{(1)},\lambda^{(2)},\dots) \ | \ \lambda^{(i)} \vdash \mu_i \ (i \ge 1)\},
		$$ 
		so that, by Proposition \ref{prop:SumRC}, 
		$$
		F_\mu^k= \prod_{i \ge 1}  \sum_{\lambda^{(i)} \vdash \mu_i} \RC(\lambda^{(i)}) 
		= \prod_{i \ge 1} \Cat_{\mu_i}.
		$$
	\end{proof}

	\noindent
	\emph{Remark}:
	For the double covering $\widetilde{\S}(n)$ of the symmetric group,
	a result similar to Theorem \ref{cor:complete} was obtained by Tysse and Wang \cite{TW}.
	They deal with $e_k(M_1^2,\dots,M_n^2)$, where the $M_i$ are elements of the spin group
	algebra of $\widetilde{\S}(n)$ called  
	\emph{odd Jucys-Murphy elements}.  \qed

	Corollary \ref{cor:complete} was first obtained by Murray \cite[Corollary 6.4]{Murray}
	in the framework of the Farahat-Higman algebra, and independently rediscovered by
	the second author \cite{Novak:BC} via Collins' work \cite{Collins} on unitary matrix integrals.
	The proof given here is different from either of these, and is completely combinatorial.  In fact, 
	when combined with the $1/N$ expansion \eqref{eqn:fundamental} obtained in \cite{Novak:BC},
	it yields an elementary and transparent proof of the first order asymptotics \eqref{eqn:Asymptotics} of the basic inner products
	$\langle u_{11} \dots u_{nn},u_{1\pi(1)} \dots u_{n\pi(n)} \rangle_N$. 
		
	\begin{corollary}
		\label{cor:FirstOrder}
		Let $N \geq n$ be positive integers.  For any $\pi \in C_\mu(n)$, we have
			
			\begin{equation*}
				(-1)^{|\mu|} N^{n+|\mu|} \langle u_{11} \dots u_{nn}, u_{\pi(1)} \dots u_{n\pi(n)} \rangle = 
				\prod_{i=1}^{\ell(\mu)} \Cat_{\mu_i} + O\bigg{(} \frac{1}{N^2} \bigg{)}
			\end{equation*}
			
		\noindent
		as $N \rightarrow \infty$.
	\end{corollary} 

	\begin{proof}
		The fundamental relation \eqref{eqn:fundamental} between the basic inner products and the class expansion problem
		reads
		
			\begin{equation*}
				\langle u_{11} \dots u_{nn}, u_{1\pi(1)} \dots u_{n\pi(n)} \rangle_N = \frac{1}{N^n} \sum_{k=0}^{\infty} (-1)^k
				\frac{F_\mu^k(N)}{N^k}.
			\end{equation*}
			
		\noindent
		Since $F_\mu^k(n)$ is non-zero only if $k=|\mu|+2g$ for some integer $g \geq 0$, this becomes
		
			\begin{equation*}
				\begin{split}
					(-1)^{|\mu|}N^{n+|\mu|}\langle u_{11} \dots u_{nn}, u_{1\pi(1)} \dots u_{n\pi(n)} \rangle_N
					&= \sum_{g=0}^\infty \frac{F_\mu^{|\mu|+2g}(n)}{N^{2g}} \\
					&= F_\mu^{|\mu|} +  \sum_{g=1}^\infty \frac{F_\mu^{|\mu|+2g}(n)}{N^{2g}} \\
					&= F_\mu^{|\mu|} + \frac{\text{const}(n,\mu)}{N^2},
				\end{split}
			\end{equation*}
			
		\noindent
		where $\text{const}(n,\mu)$ is a number depending only on $n$ and $\mu$.  The result now 
		follows from Corollary \eqref{cor:complete}.
		
	\end{proof}
	
	\subsection{Proof of Theorem \ref{thm:monomial}}
	In this subsection we give the proof of Theorem \ref{thm:monomial}. 
	
	\subsubsection{Basic lemmas}
	Define the \emph{support} of a permutation $\sigma$ to be the number of points on which it acts non-trivially:  
$$
\supp(\sigma)= \{i \ | \ \sigma(i) \not=i \}.
$$
If the reduced cycle-type of $\sigma$ is $\mu$, then $|\supp(\sigma)|=\wt(\mu)$.

\begin{lemma} \label{lem:WeightPM1}
Given a permutation $\pi$ and a transposition $(s\ t)$,  
let $\Pi=\pi(s\ t)$.
Suppose that $\Lambda=(\Lambda_1,\Lambda_2,\dots)$ 
and $\lambda=(\lambda_1,\lambda_2,\dots)$ are the reduced cycle-types of $\Pi$ and $\pi$, respectively.
Then we have $|\Lambda|=|\lambda|\pm 1$.
Furthermore, if $|\Lambda|=|\lambda|+1$, then $\supp(\Pi)=\supp(\pi) \cup \{s,t\}$,
and $s,t$ belong to the same cycle of $\Pi$.
\end{lemma}

\begin{proof}
Given a permutation $\pi$ and a transposition $(s\ t)$,
the following four cases may occur:
(i) $|\supp(\pi) \cap \{s,t\}|=0$;
(ii) $|\supp(\pi) \cap \{s,t\}| = 1$;
(iii) $s,t \in \supp(\pi)$, and $s,t$ belong to different cycles of $\pi$;
(iv)  $s,t \in \supp(\pi)$, and $s,t$ belong to the same cycle of $\pi$.

For the case (i), we obtain $\Lambda=\lambda \cup (1)$ immediately.
In the case (ii), we may suppose $\supp(\pi) \cap \{s,t\} = \{s\}$.
Then $\pi$ has a cycle
$( \dots, s, \pi(s), \dots)$, and $\Pi$ has the cycle $(\dots, s, t, \pi(s), \dots)$.
Therefore $\Lambda$ has a part equal to $\lambda_j+1$.
In the case (iii), $\pi$ has two cycles of the forms 
$(\dots, \pi^{-1}(s), s, \pi(s), \dots)$ and 
$(\dots\ \pi^{-1}(t), t, \pi(t), \dots)$.
Therefore $\Pi$ has the combined cycle 
$(\dots, \pi^{-1}(s), s, \pi(t), \dots\ \pi^{-1}(t), t, \pi(s), \dots)$.
Thus, a certain part $\Lambda_k$ of $\Lambda$ equals  $\lambda_i+\lambda_j+1$
for some $1 \le i<j \le \ell(\lambda)$.
In the case (iv),
$\pi$ has a cycle of the form 
$$(\dots, \pi^{-1}(s), s, \pi(s), \dots, \pi^{-1}(t), t, \pi(t), \dots),$$
and so $\Pi$ has divided cycles $(\dots, \pi^{-1}(s), s,\pi(t), \dots)$ and
$(\pi(s)\ \dots \ \pi^{-1}(t)\ t)$.
Thus, there are $\Lambda_j$ and $\Lambda_k$ equal to $r-1$ and $\lambda_i-r$ 
for some $\lambda_i$ and $r \ge 1$.

For the case (iv), $\Lambda$ and $\lambda$ satisfy the identity $|\Lambda|=|\lambda|-1$.
For other cases (i),(ii), and (iii), we have 
$|\Lambda|=|\lambda|+1$.
The rest of the claims are seen above.
\end{proof}

\begin{corollary} \label{cor:MinimalParity}
Let $\sigma$ be a permutation of reduced cycle-type $\lambda$.
Suppose that $\sigma$ factors as $(s_1,t_1) \cdots (s_p,t_p)$,
where $s_i<t_i \ (1 \le i \le p)$. Then $|\lambda| \le p$ and 
$|\lambda| \equiv p \pmod{2}$. 
\end{corollary}

If $\sigma$ is a permutation of reduced cycle-type $\lambda \vdash r$, and
if $\sigma$ may be factored into $r$ transpositions
\begin{equation} \label{eq:minif}
\sigma=(s_1,t_1) \cdots (s_r,t_r),
\end{equation}
then we say that  \eqref{eq:minif} is  a \emph{minimal} factorization of $\sigma$.

\begin{lemma} \label{lem:Support}
Let $\lambda \vdash r$ and let $\sigma$ be a permutation 
of reduced cycle-type $\lambda$.
Suppose that
$\sigma$ factors as  $(s_1,t_1)(s_2,t_2) \cdots (s_r,t_r)$,
where $s_i <t_i \ (1 \le i \le r)$ and $2 \le t_1 \le \cdots \le t_r$.
Then
$\supp(\sigma)=\{s_1,t_1,s_2,t_2,\dots,s_r,t_r\}$.
Furthermore, for each $i$, the letters $s_i,t_i$ belong to the same cycle of $\sigma$.
\end{lemma}

\begin{proof}
For each $1 \le i \le r$, define $\sigma_i=(s_1,t_1) \cdots(s_i,t_i)$.
It follows by Lemma \ref{lem:WeightPM1} that
the size of the reduced cycle-type of $\sigma_i$ must be $i$,
and that $\supp(\sigma_{i})= \supp(\sigma_{i-1}) \cup \{s_i,t_i\}$. 
In addition, $s_i,t_i$ belong to the same cycle of $\sigma_i$, 
and therefore to  the one of $\sigma$.
\end{proof}

\begin{lemma} \label{lem:DecompositionLemma}
Let $\tau^{(1)}$ and $\tau^{(2)}$ be permutations such that
$i<j$ for all $i \in \supp(\tau^{(1)})$ and $j \in \supp(\tau^{(2)})$.
Suppose that the reduced cycle-types of $\tau^{(1)}$ and $\tau^{(2)}$ have
weights $r_1$ and $r_2$, respectively.
Also, suppose that $\sigma:=\tau^{(1)}\tau^{(2)}$ may be expressed as
$\sigma=(s_1,t_1) \cdots(s_r,t_r)$, where
$r=r_1+r_2$, $s_i <t_i \ (1 \le i \le r)$, and $2 \le t_1 \le \cdots \le t_r$.
Then, 
$$
\tau^{(1)}=(s_1,t_1)\cdots(s_{r_1},t_{r_1}), \qquad
\tau^{(2)}=(s_{r_1+1},t_{r_1+1})\cdots(s_{r},t_{r}).
$$
\end{lemma}

\begin{proof}
By Lemma \ref{lem:Support}, we see 
$\supp(\tau^{(1)} ) \sqcup \supp(\tau^{(2)})
=\supp(\sigma)=\{s_1,t_1,\dots,s_r,t_r\}$.
Since $t_i$ are not decreasing, there exists an integer 
$p$ such that
$t_1,\dots,t_{p} \in \supp(\tau^{(1)})$ and $t_{p+1},\dots,t_{r} \in \supp(\tau^{(2)})$.
Furthermore, applying Lemma \ref{lem:Support} again,
we see that $s_i,t_i$ belong to the same cycle of $\sigma$, and so that
$\supp(\tau^{(1)}) = \{s_1,t_1,\dots,s_p,t_p\}$
and $\supp(\tau^{(2)}) = \{s_{p+1},t_{p+1},\dots,s_r,t_r\}$.
In particular, for any $i \in \{s_1,t_1,\dots,s_p,t_p\}$ and 
$j \in \{s_{p+1},t_{p+1},\dots,s_r,t_r\}$,
we have 
$\tau^{(1)}(i)=\sigma(i)$ and $\tau^{(2)}(j)=\sigma(j)$.

Let $\rho^{(1)}=(s_1,t_1)\cdots(s_p,t_p)$ and $\rho^{(2)}=(s_{p+1},t_{p+1})\cdots(s_r,t_r)$.
Since $\sigma =\rho^{(1)} \rho^{(2)}$
we have $\{s_1,t_1,\dots,s_p,t_p\} = \supp(\rho^{(1)})$
and $\{s_{p+1},t_{p+1},\dots,s_r,t_r\} = \supp(\rho^{(2)})$. 
Therefore for any $i \in \{s_1,t_1,\dots,s_p,t_p\}$ and 
$j \in \{s_{p+1},t_{p+1},\dots,s_r,t_r\}$, we have
$\rho^{(1)}(i)=\sigma(i)$ and $\rho^{(2)}(j)=\sigma(j)$.
This means $\tau^{(1)}=\rho^{(1)}$ and $\tau^{(2)}=\rho^{(2)}$.
In particular, the sizes of the reduced cycle-type of $\rho^{(1)}$ and 
$\rho^{(2)}$ are $r_1$ and $r_2$, respectively.
By definition of $\rho^{(i)}$ and Corollary \ref{cor:MinimalParity}, 
we have $r_1 \le p$ and $r_2 \le r-p$.
But $r=r_1+r_2$ so that $p=r_1$.
Therefore $\tau^{(1)}=\rho^{(1)}=(s_1,t_1)\cdots(s_{r_1},t_{r_1})$.
The desired expression for $\tau^{(2)}$ also follows.
\end{proof}

\subsubsection{Expression for cycles}

Let $a,r$ be non-negative integers.
Define the set $\E(a;r)$ by
$$
\E(a;r)=\{(i_1,\dots,i_r) \in \mathbb{Z}^r \ | \ 
i_1 \le \cdots \le i_r, \quad i_p \ge a+p \ (1 \le p \le r-1), \quad i_r=a+r\}
$$
for $r \ge 1$ and let $\E(a;0)=\emptyset$.
This extends the above definition of $\E(r)=\E(0;r)$, and
the mapping $(i_1,\dots,i_r) \mapsto (a+i_1,\dots,a+i_r)$ gives
a bijection from $\E(r)$ to $\E(a;r)$. 
Put 
\begin{align*}
\E_0(a;r)=& \{(i_1,\dots,i_r) \in \E(a;r) \ | \ i_p>a+p \ (1 \le p \le r-1)\}, \\
\E_1(a;r)=& \{(i_1,\dots,i_r) \in \E(a;r) \ | \ i_1=a+1, \ i_p>a+p \ (2 \le p \le r-1)\}, \\
\vdots & \\
\E_q(a;r)=& \{(i_1,\dots,i_r) \in \E(a;r) \ | \ i_q=a+q, 
\ i_p>a+p \ (q+1 \le p \le r-1)\}, \\
\vdots & \\
\E_{r-1}(a;r)=& \{(i_1,\dots,i_r) \in \E(a;r) \ | \ i_{r-1}=a+r-1 
\}.
\end{align*}
Then we obtain the decomposition $\E(a;r)= \bigsqcup_{q=0}^{r-1} \E_q(a;r)$.
For each $(i_1,\dots,i_{r}) \in \E_q(a;r)$ with $0 \le q \le r-2$,
we have $i_{r-1}=a+r$.
Therefore,
for each $0 \le q \le r-1$, the mapping
$$
(i_1,\dots,i_q,i_{q+1},\dots,i_r) \mapsto
((i_1,\dots,i_q), (i_{q+1},\dots,i_{r-1}))
$$
gives a bijection from $E_q(a;r)$ to $\E(a;q) \times \E(a+q+1;r-1-q)$.
Here when either $q=0$ or $q=r-1$, 
we regard the set $\E(a;q) \times \E(a+q+1;r-1-q)$
as $\E(a+1;r-1)$ or $\E(a;r-1)$, respectively.
Thus, we obtain a natural identification
\begin{equation} \label{eq:Eidentification}
\begin{array}{ll}
\E(a;r) 
& = \E_0(a;r) \sqcup \bigg{(} \bigsqcup_{q=1}^{r-2} \E_q(a;r) \bigg{)} \sqcup \E_{r-1}(a;r) \\
& \cong  \E(a+1;r-1) \sqcup 
\bigg{(} \bigsqcup_{q=1}^{r-2} (\E(a;q) \times \E(a+q+1;r-1-q)) \bigg{)} \sqcup \E(a;r-1).
\end{array}
\end{equation}
In particular, $|\E(r)| = |\E(r-1)|+
\sum_{q=1}^{r-2} |\E(q)| |\E(r-1-q)|+|\E(r-1)|$ for $r \ge 2$.
Comparing this equation with the Catalan recurrence,
we have $|\E(a;r)|=|\E(r)|=\Cat_r$ for all $r \ge 1$.

For two positive integers $a,r$, we define the cycle $\xi(a;r)$ of length $r+1$
by
$$
\xi(a;r)=(a,a+1,\dots,a+r).
$$
For convenience, we take $\xi(a;0)$ to be the identity permutation.
The  following proposition is the key to our proof of Theorem \ref{thm:monomial}.

\begin{proposition} \label{prop:CatalanTrans}
Let $t_1,\dots,t_r$ be positive integers satisfying $2 \le t_1 \le \cdots \le t_r$.
The cycle $\xi(a;r)$ may be expressed as a product of $r$ transpositions
\begin{equation} \label{eq:ShortestTrans}
\xi(a;r)= (s_1, t_1)(s_2, t_2) \cdots (s_r, t_r), \qquad s_i<t_i \ (1 \le i \le r)
\end{equation}
if and only if 
\begin{equation} \label{eq:AssumptionT}
(t_1,\dots,t_r) \in \E(a;r).
\end{equation}
Furthermore, for each $(t_1,\dots,t_r) \in \E(a;r)$,
the expression \eqref{eq:ShortestTrans} of $\xi(a;r)$ is unique.
\end{proposition}

\begin{example}
Consider the cycle $\xi(1;9)=(1,2,\dots,10)$
and three sequences 
$$
(3,5,5,5,8,8,8,9,10), \quad (3,4,4,7,7,9,9,10,10), \quad
(9,9,9,9,10,10,10,10,10) 
$$
in $\E(1;9)$.
The corresponding expressions of $\xi(1;9)$ are given as follows:
\begin{align*}
& (2,3)(4,5)(3,5)(1,5)(7,8)(6,8)(5,8)(8,9)(9,10), \\
& (2,3)(3,4)(1,4)(6,7)(5,7)(8,9)(7,9)(9,10)(4,10),\\
& (8,9)(7,9)(6,9)(5,9)(9,10)(4,10)(3,10)(2,10)(1,10).
\end{align*}
\end{example}

\begin{proof}[Proof of Proposition \ref{prop:CatalanTrans}]
We proceed by induction on $r$.
When $r=1$, since $\xi(a;1)=(a,a+1)$,
and since $\E(a;1)$ consists of a sequence $(a+1)$ of length $1$,
our claims are trivial.
Let $r >1$ and suppose that for cycles of length $<r+1$, all claims in the theorem  hold true.

(i) First,
we suppose that the cycle $\xi(a;r)$ is given by the form 
\eqref{eq:ShortestTrans}.
Then we have $t_r=a+r$ because $t_r$ is the maximum among $\supp(\xi(a;r))$,
where $\supp(\xi(a;r))=\{s_1,t_1,\dots,s_r,t_r\}$ by Lemma \ref{lem:Support}.
If we write as $s_r=a+q$ with $0 \le q \le r-1$, we have
$$
(s_1,t_1)\cdots(s_{r-1},t_{r-1})= (a,a+1,\dots,a+q)(a+q+1,a+q+2,\dots,a+r).
$$
By Lemma \ref{lem:DecompositionLemma}, we see that
\begin{equation} \label{eq:twocyclesDecom}
\begin{split}
(s_1,t_1)\cdots(s_{q},t_{q})=& (a,a+1,\dots,a+q), \\
(s_{q+1},t_{q+1}) \cdots (s_{r-1},t_{r-1} )=& (a+q+1,a+q+2,\dots,a+r).
\end{split}
\end{equation}
By the induction hypothesis for cycles of length $q+1$ and of length $r-q$,
we have $(t_1,\dots,t_q) \in \E(a;q)$ and $(t_{q+1},\dots,t_{r-1}) \in \E(a+q+1;r-1-q)$.
This fact and Equation \eqref{eq:Eidentification}
imply $(t_1,\dots,t_q,t_{q+1},\dots,t_{r-1},t_r) \in \E_q(a;r) \subset \E(a;r)$.

(ii) Next, we suppose $(t_1,\dots,t_r) \in \E(a;r)$.
According to the decomposition $\E(a;r)=\bigsqcup_{q=0}^{r-1} 
\E_q(a;r)$,
there exists a unique number $q$ such that 
$0 \le q \le r-1$ and $(t_1,\dots,t_r) \in \E_q(a;r)$, and then
$(t_1,\dots,t_q) \in \E(a;q)$ and $(t_{q+1},\dots,t_{r-1}) \in \E(a+q+1;r-1-q)$.
By the induction assumption, there exist  sequences 
$(s_1,s_2,\dots,s_{q})$ and $(s_{q+1},\dots,s_{r-1})$ satisfying
\eqref{eq:twocyclesDecom}.
Therefore we obtain the expression
$$
\xi(a;r)=(s_1,t_1) \cdots (s_{q},t_q)(s_{q+1},t_{q+1})\cdots(s_{r-1},t_{r-1})(a+q,a+r),
$$
as required.

(iii) It remains to prove the uniqueness of the expression \eqref{eq:ShortestTrans}.
Assume that 
the cycle $\xi(a;r)$ has two expressions
$$
(s_1, t_1)(s_2, t_2) \cdots (s_r, t_r) \qquad \text{and}
\qquad (s_1',t_1)(s_2',t_2) \cdots (s_r',t_r), 
$$
where $s_i,s_i'<t_i \ (1 \le i \le r)$.
Write as $s_r=a+q$ and $s_r'=a+q'$.
As we saw in the part (i), 
the sequence $(t_1,\dots,t_r)$ belongs to $\E_q(a;r) \cap \E_{q'}(a;r)$.
But, since $\E_q(a;r) \cap \E_{q'}(a;r)= \emptyset$ if $q \not=q'$, 
we have $q=q'$ so that $s_r=s_r'$.
Now, as like  \eqref{eq:twocyclesDecom}, 
we have $(t_1,\dots,t_q) \in \E(a;q)$ and $(t_{q+1},\dots,t_{r-1}) \in \E(a+q+1;r-1-q)$,
and
\begin{equation*}
\begin{split}
(s_1,t_1)\cdots(s_{q},t_{q})=(s_1',t_1)\cdots(s_{q}',t_{q})=& (a,a+1,\dots,a+q), \\
(s_{q+1},t_{q+1}) \cdots (s_{r-1},t_{r-1} )=(s_{q+1}',t_{q+1}) \cdots (s_{r-1}',t_{r-1} )
=& (a+q+1,a+q+2,\dots,a+r).
\end{split}
\end{equation*}
By the induction assumption, we obtain
$s_1=s_1',\dots, s_q=s_q', s_{q+1}=s'_{q+1},\dots, s_{r-1}=s_{r-1}'$.
\end{proof}

\subsubsection{Proof of Theorem \ref{thm:monomial}} \label{subsection:ProofEnd}

Recall the definition of the Jucys-Murphy elements: $J_t=\sum_{1 \le s <t} (s,t)$.
For a permutation $\sigma \in S(n)$ and  a polynomial $f$ in $n$ variables, 
denote by $[\sigma] f(\Xi_n)$
the multiplicity of $\sigma$ in $f(J_1,\dots,J_n)$:
$$
f(\Xi_n) = \sum_{\sigma \in S(n)} \bigg{(} [\sigma] f(\Xi_n)\bigg{)} \sigma \in \C[S_n].
$$
For a partition $\mu$ with size $k$ and length $l$, we define the canonical permutation $\sigma_\mu$ 
of reduced cycle-type $\mu$ by
\begin{align*}
\sigma_\mu=& (1,2,\dots,\mu_1+1)(\mu_1+2,\dots, \mu_1+\mu_2+2)\cdots
(\mu_1+\cdots+\mu_{l-1}+l, \dots, k+l) \\
=& \xi(1;\mu_1) \xi(\mu_1+2;\mu_2) \cdots \xi(\mu_1+\cdots+\mu_{l-1}+l;\mu_l).
\end{align*}

\begin{proposition} \label{prop:DecomposeCatalan}
Let $\mu$ be a partition of $k$
and let $(t_1,\dots,t_k)$ be a sequence of positive integers such that 
$2 \le t_1 \le \cdots \le t_k$. 
Then  $[\sigma_\mu] J_{t_1} \cdots J_{t_k} =1$ if 
$(t_1,\dots,t_k)$ satisfies
\begin{equation} \label{eq:TiCatalans}
(t_{\mu_1+\cdots+\mu_{i-1}+1},\dots, t_{\mu_1+\cdots+\mu_{i-1}+\mu_i}) \in
\E(\mu_1+\cdots+\mu_{i-1}+i;\mu_i)
\end{equation}
for all $1 \le i \le \ell(\mu)$,
and $[\sigma_\mu] J_{t_1} \cdots J_{t_k} =0$ otherwise.
\end{proposition}

\begin{proof}
The value $[\sigma_\mu] J_{t_1} \cdots J_{t_k}$ is the number of
sequences $(s_1,\dots,s_k)$ satisfying
$$
\sigma_\mu = 
\prod_{i=1}^{\ell(\mu)} \xi(\mu_1+\dots+\mu_{i-1}+i;\mu_i) = (s_1,t_1) \cdots(s_k,t_k).
$$
By Lemma \ref{lem:DecompositionLemma}, it equals 
the number of sequences $(s_1,\dots,s_k)$ satisfying
$$
\xi(\mu_1+\dots+\mu_{i-1}+i;\mu_i)
=(s_{\mu_1+\cdots+\mu_{i-1}+1},t_{\mu_1+\cdots+\mu_{i-1}+1}) \cdots 
(s_{\mu_1+\cdots+\mu_{i-1}+\mu_i},t_{\mu_1+\cdots+\mu_{i-1}+\mu_i}) 
$$
for all $1 \le i \le \ell(\mu)$.
It follows by Proposition \ref{prop:CatalanTrans} that
$[\sigma_\mu] J_{t_1}\cdots J_{t_k}$ equals to $1$ if
\eqref{eq:TiCatalans} holds true for all $i$, and 
to $0$ otherwise.
\end{proof}

\begin{example}
Let $2 \le t_1 \le \cdots \le t_6$ and consider 
$\sigma_{(3,2,1)}=(1,2,3,4)(5,6,7)(8,9)$.
Suppose $[\sigma_{(3,2,1)}] J_{t_1} \cdots J_{t_{6}} =1$.
Then, Proposition \ref{prop:DecomposeCatalan} claims 
$$
(t_1,t_2,t_3) \in \E(1;3), \qquad (t_4,t_5) \in \E(5;2), \qquad
(t_6) \in \E(8;1).
$$
Therefore, $(t_1,t_2) \in \{(2,3), (2,4),(3,3), (3,4), (4,4)\}$, $t_3=4$, 
$t_4\in \{6,7\}$, $t_5=7$, and $t_6=9$.
\end{example}

As defined above,
a weakly increasing sequence $(t_1,\dots,t_r)$ is of type $\lambda \vdash r$ with $\ell(\lambda)=l$
if there exists a permutation $(\alpha_1,\dots,\alpha_l)$ of 
$(\lambda_1, \dots, \lambda_l)$
such that
$$
t_{1}= t_{2}  = \cdots = t_{\alpha_1} <
t_{\alpha_1+1}= t_{\alpha_1+2}  = \cdots = 
t_{\alpha_1+\alpha_2}<
t_{\alpha_1+\alpha_2+1}= \cdots.
$$  
The monomial symmetric polynomial $m_{\lambda}(\Xi_n), \ \lambda \vdash k,$
is written as
$$
m_{\lambda}(\Xi_n)= \sum_{\begin{subarray}{c}
2 \le t_1 \le \cdots \le t_k \le n \\
(t_1,\dots,t_k): \text{type $\lambda$} 
\end{subarray}} 
J_{t_1} J_{t_2} \cdots J_{t_k}
= \sum_{\begin{subarray}{c}
2 \le t_1 \le \cdots \le t_k \le n \\
(t_1,\dots,t_k): \text{type $\lambda$} 
\end{subarray}} 
\sum_{s_1=1}^{t_1-1} \cdots \sum_{s_k=1}^{t_k-1} (s_1,t_1) \cdots (s_k,t_k).
$$

Let $\mu$ be a partition of $k$.
We now evaluate 
the coefficient $L^\lambda_\mu(n)$ of $\cc_\mu(n)$ in $m_\lambda(J_1,\dots,J_n)$,
which equals $L_\mu^\lambda(n)=[\sigma_\mu] m_\lambda(J_1,\dots,J_n)$.
By the assumption $n \ge k+\ell(\mu)$, the permutation $\sigma_\mu$ lives in $S(n)$.
It follows by Proposition \ref{prop:DecomposeCatalan} that
$L^\lambda_\mu(n)$ is the number of weakly increasing 
sequences $(t_1,\dots,t_k)$ of type $\lambda$, satisfying \eqref{eq:TiCatalans}
for all $1 \le i \le \ell(\mu)$.
If $(t_1,\dots,t_k)$ is such a sequence and if we 
let $\lambda^{(i)} \vdash \mu_i$ being the type of $(t_{\mu_1+\cdots+\mu_{i-1}+1},
\dots,t_{\mu_1+\cdots+\mu_{i-1}+\mu_i})$, then
$\lambda$ must agree with $\lambda^{(1)} \cup \lambda^{(2)} \cup \cdots$
so that $(\lambda^{(1)},\lambda^{(2)},\dots) \in \mathfrak{R}(\lambda,\mu)$.
Thus, $L^\lambda_\mu(n)$ coincides with
\begin{align*}
& \sum_{(\lambda^{(1)},\lambda^{(2)},\dots) \in \mathfrak{R}(\lambda,\mu)}
\prod_{i =1}^{\ell(\mu)} (\text{the number of sequences 
in $\E(\mu_1+\cdots+\mu_{i-1}+i;\mu_i)$
of type $\lambda^{(i)}$}) \\
=& \sum_{(\lambda^{(1)},\lambda^{(2)},\dots) \in \mathfrak{R}(\lambda,\mu)}
\prod_{i =1}^{\ell(\mu)} \RC(\lambda^{(i)}).
\end{align*}
This completes the proof of Theorem \ref{thm:monomial}.

\section{Class coefficients and connection coefficients}
\label{sec:Connection}

In this section, we pursue a certain analogy between the class coefficients $G_\mu(f,n)$ arising in the expansion of symmetric 
functions in Jucys-Murphy elements, and the \emph{connection coefficients} $A^{\alpha\beta}_{\mu}(n)$ of the class
algebra $\Z(n)$.  The latter are by definition the structure constants of $\Z(n)$, i.e.

	\begin{equation*}
		\cc_\alpha(n) \cc_\beta(n) = \sum_{\mu} A^{\alpha\beta}_{\mu}(n) \cc_\mu(n).
	\end{equation*}
	
\noindent
Like the class coefficient $L_\mu^\lambda(n)$, the connection coefficient $A_\mu^{\alpha\beta}(n)$ has an immediate combinatorial 
interpretation: it is equal to the number of factorizations

	\begin{equation*}
		\pi = \sigma \rho
	\end{equation*}
	
\noindent
of a fixed representative $\pi$ of $C_\mu(n)$ into a permutation $\sigma$ of reduced cycle type $\alpha$ and a permutation 
$\rho$ of reduced cycle type $\beta$.  The following properties of connection coefficients are well known:

	\begin{enumerate}
		
		\item
		$A_\mu^{\alpha\beta}(n)$ is a polynomial function of $n$;
			
			\medskip
			
		\item
		$A_\mu^{\alpha\beta}(n)$ vanishes unless $|\mu|+2g = |\alpha|+|\beta|$ for some integer $g \geq 0$;
		
			\medskip
			
		\item
		In the ``top'' case $|\mu|=|\alpha|+|\beta|$, the connection coefficient $A_\mu^{\alpha\beta}=A_\mu^{\alpha\beta}(n)$
		is independent of $n$, and vanishes unless $\alpha \cup \beta$ is a refinement of $\mu$.
		
	\end{enumerate}

Property $(1)$ above is a classical result due to Farahat and Higman \cite{FH}, see also \cite{IK}.
In the previous section, we proved analogues of properties $(2)$ and $(3)$ for the class coefficients $L_\mu^\lambda(n)$.
In this section, we will prove the analogue of $(1)$, namely that for any fixed symmetric function $f$ and partition $\mu$,
the class coefficient $G_\mu(f,n)$ is a polynomial function of $n$.  By combining this fact with the spectral properties of Jucys-Murphy 
elements in irreducible symmetric group representations, we recover and generalize a recent result of Stanley \cite{Stanley:content}
and Olshanski \cite{Olshanski} on the polynomiality
of certain statistics on Plancherel-random partitions.

Macdonald \cite{Macdonald}, see also \cite{GJ:top}, used Lagrange inversion to construct a basis $\{g_\mu\}$ of the algebra
of symmetric functions which encodes the top class coefficients $A_\mu^{\alpha\beta}$:

	\begin{equation}
		\label{eqn:Macdonald}
		g_\alpha g_\beta = \sum_{|\mu|=|\alpha|+|\beta|} A_\mu^{\alpha\beta} g_\mu.
	\end{equation}
	
\noindent
We conclude this section by obtaining an analogue of Macdonald's result, which realizes the top class coefficients $L_\mu^\lambda$
intrinsically as part of the algebraic structure of $\Lambda$, without any reference to Jucys-Murphy elements.  This computation 
constitutes a new change of basis formula in the algebra of symmetric functions.

	\subsection{Polynomiality}
			
		\begin{theorem}
		\label{thm:Polynomial}
			Fix a symmetric function $f \in \Lambda$ and a partition $\mu$.  The class coefficient $G_\mu(f,n)$ is a polynomial
			function of $n$.
		\end{theorem}
			
		\begin{proof}
		Since $\Lambda=\C[e_1,e_2,\dots]$ there exists a polynomial, say $p_f,$ such that
			$$f=p_f(e_{i_1},\dots,e_{i_k})$$
		for some elementary symmetric functions $e_{i_1},\dots,e_{i_k}.$
		By Proposition \ref{prop:Jucys} we have
			$$f(\Xi_n) = p_f\bigg{(}
			\sum_{|\mu|=i_1} \cc_{\mu}(n), \dots, \sum_{|\mu|=i_k} \cc_{\mu}(n)
			\bigg{)},$$
		and the result now follows from the polynomiality of the connection coefficients $A_\mu^{\alpha\beta}(n)$.
		\end{proof}

		From the isotypic decomposition 
		
			\begin{equation*}
				\C[\S(n)] = \bigoplus_{\lambda \vdash n} (\dim \lambda) \mathsf{V}^\lambda,
			\end{equation*}
			
		\noindent
		where $\dim \lambda$ denotes the dimension of the irreducible representation $\mathsf{V}^\lambda$,
		one obtains \emph{Burnside's identity}:
		
			\begin{equation*}
				\sum_{\lambda \vdash n} (\dim \lambda)^2 = n!.
			\end{equation*}
			
		\noindent
		Burnside's identity implies that the function
		
			\begin{equation*}
				\lambda \mapsto \frac{(\dim \lambda)^2}{n!}
			\end{equation*}
			
		\noindent
		defines a probability measure, known as the \emph{Plancherel measure}, on the sample space $\mathcal{Y}_n=\{\lambda \vdash n\}$.  
		Given a function $f:\mathcal{Y} \rightarrow \C$ defined on the set of partitions, let
		
			\begin{equation*}
				\langle f \rangle_n = \sum_{\lambda \vdash n} f(\lambda) \frac{(\dim \lambda)^2}{n!}
			\end{equation*}
			
		\noindent
		denote its expected value with respect to the Plancherel measure on $\mathcal{Y}_n$.  The following polynomiality property of 
		Plancherel averages was proved by Stanley \cite{Stanley:content} and Olshanski \cite{Olshanski}, by different methods.
		As an application of Theorem \ref{thm:Polynomial}, we give a third proof.
		
		\begin{theorem}
			Let $f \in \Lambda$ be a symmetric function.  The Plancherel expectation 
			$\langle f(A_\lambda) \rangle_n$, where $A_\lambda$ is the content alphabet of $\lambda$, is 
			a polynomial function of $n$.
		\end{theorem}
		
		\begin{proof}
			We will prove the following more general fact.  Let $\chi^\lambda_\mu$ denote the trace of any
			representative of the conjugacy class $C_\mu(n)$ in the irreducible representation $\mathsf{V}^\lambda$.
			Then, the sum
			
				\begin{equation*}
					\sum_{\lambda \vdash n} f(A_\lambda) \chi^\lambda_\mu \frac{\dim \lambda}{n!}
				\end{equation*}
				
			\noindent
			is a polynomial function of $n$.  When $\mu$ is the empty partition, this reduces to the statement
			of the theorem.
			
			We prove the more general assertion as follows.  Put
			
				\begin{equation*}
					\chi^\lambda = \sum_\mu \chi^\lambda_\mu \cc_\mu(n).
				\end{equation*}
				
			\noindent
			Then it is a basic fact that $\{\chi^\lambda : \lambda \vdash n\}$ constitutes a basis of the class algebra
			$\Z(n)$, and the coordinates of the identity are given in this basis by
			
				\begin{equation*}
					\cc_{(0)}(n) = \sum_{\lambda \vdash n} \frac{\dim \lambda}{n!}\chi^\lambda.
				\end{equation*}
				
			\noindent
			Since the central character of $f(\Xi_n)$ in $\mathsf{V}^\lambda$ is $f(A_\lambda)$, we have 
			$f(\Xi_n)\chi^\lambda = f(A_\lambda)\chi^\lambda$ in $\Z(n)$.  Thus the coordinates of $f(\Xi_n)$ relative
			to the character basis are
			
				\begin{equation*}
					f(\Xi_n) = f(\Xi_n)\cc_{(0)}(n) = \sum_{\lambda \vdash n} f(A_\lambda) \frac{\dim \lambda}{n!} \chi^\lambda,
				\end{equation*}
				
			\noindent
			so that
			
				\begin{equation*}
					G_\mu(f,n) = \sum_{\lambda \vdash n} f(A_\lambda) \frac{\dim \lambda}{n!}\chi^\lambda_\mu.
				\end{equation*}
				
			\noindent
			The result now follows from Theorem \ref{thm:Polynomial}.
		\end{proof}
		
		\subsection{An analogue of Macdonald's result for top connection coefficients}
		Macdonald \cite[Chapter I.7, Example 25]{Macdonald}, see also \cite{GJ:top}, 
		used Lagrange inversion to
		construct a basis $\{g_{\mu}\}$
		of the algebra of symmetric functions whose connection coefficients coincide with the 
		top connection coefficients $A_{\alpha\beta}^{\mu}$ of the class algebra, as in equation
		\eqref{eqn:Macdonald}.

	We now give an analogue of Macdonald's result for the top 
	class coefficients $L_{\mu}^{\lambda}:$ for each $k\geq 1$ we realize
	the matrix $(L_{\mu}^{\lambda})_{|\mu|=|\lambda|=k}$ as the transition matrix between two
	bases of the degree $k$ component of the graded algebra $\Lambda.$
	
	Since the elementary symmetric functions $e_k$ are algebraically independent and
	generate $\Lambda,$ we may define an endomorphism $\psi:\Lambda \rightarrow \Lambda$
	by $\psi(e_k)=h_k.$  This endomorphism is in fact involutive: $\psi(h_k)=e_k.$
	The image $f_{\lambda} := \psi(m_{\lambda})$ of the monomial symmetric function of 
	type $\lambda$ under $\psi$ 
	is known as the \emph{forgotten symmetric function} of type $\lambda,$ see 
	\cite[Exercise 7.9]{Stanley:EC2}.
	
	\begin{theorem}
		\label{thm:forgotten}
		Let $|\lambda|=k.$  Then
			$$(-1)^{k} f_{\lambda} = 
			\sum_{|\mu|=k} L_{\mu}^{\lambda} g_{\mu}.$$
	\end{theorem}
	
	\begin{proof}
		Let 
			$$
			u= t+\sum_{r=1}^\infty h_r t^{r+1}.
			$$
		Then $t$ can be expressed as a power series in $u$.
		Define symmetric functions $h_r^*$, $r=1,2,\dots$, via
			$$
			t=u+\sum_{r=1}^\infty h_r^* u^{r+1}.
			$$
		From the Lagrange inversion formula,
		the symmetric functions $h_r^*$ are explicitly given by
		\begin{equation*}
		h_r^*= (-1)^r \sum_{\lambda \vdash r} \RC(\lambda) e_\lambda,
		\end{equation*}
		where $e_\lambda=e_{\lambda_1}e_{\lambda_2} \cdots$,
		see \cite[(3.6)]{Murray} and also \cite{Haiman},  \cite[Ch. I, Example 2.24]{Macdonald},
		\cite{Stanley:park}.

		Let $h_\lambda^*=h_{\lambda_1}^* h_{\lambda_2}^* \cdots$.
		Then $\{h_\lambda^*\}$ is a basis of $\Lambda$.
		
		Let $\mu$ be a partition of $k$. We will now prove that
			\begin{equation} \label{eq:TransitionL}
			(-1)^k h^*_\mu= \sum_{\lambda \vdash k } L^\lambda_\mu e_\lambda,
			\end{equation}
		and thus the matrix $(L^\lambda_\mu)_{|\mu|=|\lambda|=k}$ is 
		the transition matrix from the basis $\{(-1)^{|\lambda|} h^*_\lambda\}$ to
		the basis $\{e_\lambda\}$ in the $k$th component of $\Lambda.$
		The proof goes as follows: let $l=\ell(\mu)$. It follows from Theorem \ref{thm:monomial} that
		\begin{align*}
		h_\mu^* =& h_{\mu_1}^* h_{\mu_2}^* \cdots h_{\mu_l}^* \\
		=& (-1)^k \sum_{\lambda^{(1)} \vdash \mu_1} \sum_{\lambda^{(2)} \vdash \mu_2} 
		\cdots \sum_{\lambda^{(l)} \vdash \mu_l} 
		\RC(\lambda^{(1)}) \RC(\lambda^{(2)}) \cdots \RC(\lambda^{(l)}) 
		e_{\lambda^{(1)}\cup \lambda^{(2)} \cup \cdots \cup \lambda^{(l)}} \\
		=& (-1)^k \sum_{\lambda \vdash k} \sum_{(\lambda^{(1)},\lambda^{(2)},\dots,\lambda^{(l)})
		\in \R(\lambda,\mu)} \RC(\lambda^{(1)}) \RC(\lambda^{(2)}) \cdots \RC			(\lambda^{(l)}) 
			e_\lambda \\
		=& (-1)^k\sum_{\lambda \vdash k} L^\lambda_\mu e_\lambda.
		\end{align*}

		Let $\langle \cdot, \cdot \rangle$ be the scalar product on $\Lambda$ defined by
		$\langle h_\lambda, m_\mu \rangle = \delta_{\lambda \mu}$.
		With respect to this scalar product, the dual bases of $\{h_\lambda^*\}$
		and $\{e_\lambda\}$ are, respectively, Macdonald's symmetric functions
		$\{g_\lambda\}$ and the forgotten symmetric functions $\{f_\lambda\}$,
		see  \cite[Ch. I.2]{Macdonald},  \cite[Ch. I, Example 7.25]{Macdonald},  \cite{Murray}.
		Thus \eqref{eq:TransitionL} is equivalent to 
		\begin{equation*}
		(-1)^k f_\lambda= \sum_{\mu \vdash k } L^\lambda_\mu g_\mu, \qquad \lambda \vdash k.
		\end{equation*}
		
	\end{proof}
	
\section{Coefficient of a full cycle and central factorial numbers}
\label{sec:CentralFactorial}

	Given a permutation $\pi \in \S(n)$ and an integer $k \geq 0$, how many factorizations
	
		\begin{equation*}
			\pi = (s_1\ t_1) (s_2\ t_2) \dots (s_k\ t_k)
		\end{equation*}
		
	\noindent
	of $\pi$ into $k$ transpositions are there?  This is a very natural question.  Geometrically,
	the problem is to count the number of $k$-step walks from the identity to $\pi$ on the Cayley graph of $\S(n)$.
	Algebraically, this is a special case of the connection coefficient problem, since if $\pi \in C_\mu(n)$ then the desired number, call it
	$A_\mu^k(n)$, is 
	the coefficient of $\cc_\mu(n)$ in the product $\cc_{(1)}(n)\cc_{(1)}(n) \dots \cc_{(1)}(n)$ of $k$ copies of 
	the class of transpositions.  In the special case where $\pi \in C_{(n-1)}(n)$ is a full cycle, the number of 
	factorizations of $\pi$ into $n-1$ transpositions (i.e. the minimal number required) is the Cayley number,
	
		\begin{equation}
			\label{eqn:FullCycleMinimal}
			A_{(n-1)}^{n-1}(n)=\operatorname{Cay}_{n-1}=n^{n-2},
		\end{equation}
	
	\noindent
	which the reader will recognize as the number of trees on the 
	vertex set $[n-1]$.  This formula was first obtained by Hurwitz \cite{Hurwitz} as a corollary
	of his exact enumeration
	of holomorphic maps $\mathbb{P}^1 \rightarrow \mathbb{P}^1$ from the Riemann sphere to itself with one degenerate branch point.
	Various bijective proofs of this result have since been found, see \cite{GP} and references therein.  It follows that
	the number of factorizations of a fixed representative of $C_\mu(n)$ into $|\mu|$ transpositions, the minimal number
	required, is 
	
		\begin{equation}
			\label{eqn:Cayley}
			A_{\mu}^{|\mu|} = {|\mu| \choose \mu_1,\dots,\mu_{\ell(\mu)}} \prod_{i=1}^{\ell(\mu)} \operatorname{Cay}_{\mu_i}.
		\end{equation}
	
	\noindent
	The number $A_{(n-1)}^k(n)$ of factorizations of a fixed representative of the class $C_{(n-1)}(n)$ of full cycles into any number $k$
	of transpositions was determined by Jackson \cite{Jackson} using a character-theoretic argument.  The number of factorizations is 
	non-zero if and only if $k=n-1+2g$ for some integer $g \geq 0$, and Jackson's formula is
		
		\begin{equation}
			\label{eqn:Jackson}
			A_{(n-1)}^{n-1+2g} = \frac{1}{n!} \sum_{k=0}^{n-1} (-1)^k {n-1 \choose k} \bigg{(} {n \choose 2}-kn \bigg{)}^{n-1+2g}.
		\end{equation}
		
	\noindent
	Jackson's formula
	was rediscovered in the context of singularity theory by Shapiro, Shapiro and Vainshtein \cite{SSV}, and a bijective
	proof was found by Goulden \cite{Goulden}.
	
	Consider now the following restricted variant of the above question: given a permutation $\pi$ and an integer
	$k \geq 0$, how many factorizations
	
		\begin{equation*}
			\pi = (s_1\ t_1) (s_2\ t_2) \dots (s_k\ t_k)
		\end{equation*}
		
	\noindent
	of $\pi$ into $k$ transpositions are there satisfying the constraint $t_1 \leq t_2 \leq \dots \leq t_k$?  Factorizations of this
	form were called \emph{primitive} in \cite{GM}.  Algebraically,
	this is a special case of the class expansion problem: if $\pi \in C_\mu(n),$ then the required number is the 
	coefficient $F_\mu^k(n)$ of $\cc_\mu(n)$ in $h_k(\Xi_n)$, where $h_k$ is the complete homogeneous symmetric function 
	of degree $k$.  
	For the class of full cycles, the analogue of \eqref{eqn:FullCycleMinimal} was first obtained by 
	Gewurz and Merola \cite{GM}:
	
		\begin{equation}
			F_{(n-1)}^{(n-1)} = \Cat_{n-1}.
		\end{equation}
		
	\noindent
	The general solution to this problem in the minimal case $k=|\mu|$ is given by Corollary \ref{cor:complete},
	which gives the analogue of \eqref{eqn:Cayley}:
	
		\begin{equation}
			F_{\mu}^{|\mu|} = \prod_{i=1}^{\ell(\mu)} \Cat_{\mu_i}.
		\end{equation}
		
	\noindent
	Note that the order constraint has the effect of desymmetrizing the Cayley numbers to the Catalan numbers and removing
	the shuffle factor.  In this language, Theorem \ref{thm:monomial} is an exact enumeration of minimal primitive factorizations by type.
	
	In this final section, we obtain the analogue of Jackson's formula \eqref{eqn:Jackson}.  That is, we give 
	an exact enumeration of the primitive factorizations of any full
	cycle, of any length.  Like Jackson's argument, our method relies on techniques from the character theory of the symmetric groups.
	We prove that
	
		\begin{equation}
			\label{eqn:CentralFactorial}
			F_{(n-1)}^{n-1+2g}(n) = \Cat_{n-1} \cdot T(n-1+g,n-1)
		\end{equation}

	\noindent
	for any integer $g \geq 0$, where
	
		\begin{equation}
			T(a,b) = 2 \sum_{j=1}^b (-1)^{b-j}\frac{j^{2a}}{(b-j)!(b+j)!}
		\end{equation}
		
	\noindent
	is the Carlitz-Riordan \emph{central factorial number} of the second kind.
	$T(a,b)$ is equal to the number of partitions of the
	set $\{1,1',2,2',\dots,a,a'\}$ into $b$ disjoint non-empty subsets $V_1,\dots,V_b$ such that, for each 
	$1 \leq k \leq b$, if $i$ is the least integer such that either $i$ or $i'$ belongs to $V_k$ then
	$\{i,i'\} \subseteq V_k$.  Thus $T(a,b)$ may be thought of as a two-coloured Stirling number of the second kind.  
	Indeed, central factorial numbers first appeared in Carlitz and Riordan's investigation of the central 
	difference operator \cite{CR}, where they play the same role as the Stirling numbers for the usual difference operator.
	For more information and references regarding central factorial numbers, see \cite[Section 6.5]{Riordan} and \cite[Exercise 5.8]{Stanley:EC2}.
	The result presented in this section gives a new and rather surprising combinatorial interpretation of the central factorial 
	numbers as counting primitive factorizations of a full cycle.  It would be very interesting to obtain a bijective proof.
	
	\subsection{Coefficient of $\cc_{(n-1)}(n)$ in $h_k(\Xi_n)$}
	Let $q$ be a formal variable, and introduce the generating function

		\begin{equation}
			\label{eqn:GeneratingFunction}
			\Phi(q,n) = \sum_{k=0}^{\infty} h_k(\Xi_n)q^k.
		\end{equation}
		
	\noindent
	The generating series $\Phi(q,n)$ is an element of the ring $\Z(n)[[q]]$ of formal power series in $q$
	with coefficients in $\Z(n)$.  Via the hook-length formula, Jucys' result \eqref{eqn:JucysCharacter} may be re-stated as follows: for any 
	symmetric function $f \in \Lambda$, the central element $f(\Xi_n) \in \Z(n)$ has coordinates
	
		\begin{equation*}	
			f(\Xi_n) = \sum_{\lambda \vdash n} \frac{f(A_\lambda)}{H_\lambda}\chi^{\lambda}
		\end{equation*}
		
	\noindent
	with respect to the character basis of $\Z(n)$, 
	where $H_\lambda = \prod_{\Box \in \lambda} h(\Box)$ is the product of the hook-lengths over the cells of $\lambda$.
	Substituting the character expansion of $h_k(\Xi_n)$ into the generating function \eqref{eqn:GeneratingFunction},
	we obtain
	
		\begin{equation*}
			\Phi(q,n) = \sum_{k=0}^{\infty} \bigg{(}  \sum_{\lambda \vdash n} \frac{f(A_\lambda)}{H_\lambda}\chi^{\lambda} \bigg{)} q^k.
		\end{equation*}
		
	\noindent
	Changing order of summation and appealing to the generating function
	
		\begin{equation*}
			\sum_{k=0}^{\infty} h_k(x_1,x_2,x_3,\dots)q^k = \prod_{i=1}^{\infty} \frac{1}{1-x_i q}
		\end{equation*}
		
	\noindent
	of the complete symmetric functions, this becomes
	
		\begin{equation}
			\label{eqn:RationalForm}
			\Phi(q,n) = \sum_{\lambda \vdash n} \frac{\chi^\lambda}{H_\lambda \prod_{\Box \in \lambda} (1-c(\Box) q)}.
		\end{equation}
		
	Now let $\mu$ be a partition, and set
	
		\begin{equation*}
			\Phi_\mu(q,n) = \sum_{k=0}^{\infty} F_\mu^k(n) q^k.
		\end{equation*}
		
	\noindent
	Then from \eqref{eqn:RationalForm} we obtain
	
		\begin{equation*}
			\Phi_\mu(q,n) = \sum_{\lambda \vdash n} \frac{\chi^\lambda_\mu}{H_\lambda \prod_{\Box \in \lambda} (1-c(\Box) q)},
		\end{equation*}	
		
	\noindent
	from which it is clear that $\Phi_\mu(q,n)$ is a rational function of $q$, with coefficients in $\mathbb{Q}$.  Now, the 
	trace of a full cycle $\pi \in C_{(n-1)}(n)$ in an irreducible representation $\mathsf{V}^\lambda$ of $\S(n)$ is
	non-zero if and only if $\lambda$ is a ``hook'' partition:
		
		\begin{equation*}
			\chi^\lambda_{(n-1)} = \begin{cases}
				(-1)^r, \text{ if } \lambda=(n-k,1^k) \\
				0, \text{ otherwise }
			\end{cases},
		\end{equation*}
		
	\noindent
	see e.g. \cite{DG} for a proof of this classical fact.  When $\lambda=(n-k,1^k)$ is a hook, the content alphabet 
	of $\lambda$ is simply
	
		\begin{equation*}
			A_\lambda = \{1,2,\dots,n-k-1\} \sqcup \{-1,-2,\dots,-(k-1)\}.
		\end{equation*}
		
	\noindent
	Thus, summing over all hook representations, we have
	
		\begin{equation*}
			\Phi_{(n-1)}(q,n) = \sum_{k=0}^{n-1} \frac{(-1)^k}{H_{(n-k,1^k)} \prod_{j=1}^{n-k-1} (1-jq) \prod_{j=1}^{k-1} (1+jq)},
		\end{equation*}
		
	\noindent
	which is a rational function of the form
	
		\begin{equation*}
			\Phi_{(n-1)}(q,n) = \frac{a_0+a_1q + \dots + a_{n-1}q^{n-1}}{\prod_{j=1}^{n-1} (1-j^2q^2)}.
		\end{equation*}
		
	\noindent
	On the other hand, we know from Corollary \ref{cor:complete} that 
	
		\begin{equation*}
			\Phi_{(n-1)}(q,n) = \sum_{g=0}^{\infty} F_{(n-1)}^{n-1+2g}(n) q^{n-1+2g}, \quad F_{(n-1)}^{n-1}=\Cat_{n-1},
		\end{equation*}
		
	\noindent
	so that $a_0=a_1=\dots=a_{n-2}=0$ and $a_{n-1}=\Cat_{n-1}$.  We have thus proved the following.
	
	\begin{theorem}
		\label{thm:CentralFactorial}
		We have
			
			\begin{equation*}
				\Phi_{(n-1)}(q,n) = \frac{\Cat_{n-1}q^{n-1}}{(1-q^2)(1-4q^2)(1-9q^2) \dots (1-(n-1)^2q^2)}.
			\end{equation*}
	\end{theorem}
	
	Up to the factor $\Cat_{n-1}$, the rational function appearing in Theorem \ref{thm:CentralFactorial} is an ordinary 
	generating function for the central factorial numbers:
	
		\begin{equation*}
			\frac{q^{n-1}}{(1-q^2)(1-4q^2)(1-9q^2) \dots (1-(n-1)^2q^2)} = \sum_{g \geq 0} T(n-1+g,n-1)q^{n-1+g},
		\end{equation*}
		
	\noindent
	see \cite[Exercise 5.8]{Stanley:EC2}.  Thus Theorem \ref{thm:CentralFactorial} is equivalent to the identity 
	\eqref{eqn:CentralFactorial} stated above.  Finally, making the substitution $q=-1/N$ and appealing to 
	the fundamental identity \eqref{eqn:fundamental}, we see that Theorem \ref{thm:CentralFactorial} 
	implies the following exact formula \cite{Collins} for cyclic inner products in the algebra $\mathcal{A}$:
	
		\begin{equation}
			\langle u_{11}u_{22}u_{33} \dots u_{nn}, u_{12}u_{23}u_{34} \dots u_{n1} \rangle_N
			= \frac{(-1)^{n-1}\Cat_{n-1}}{N(N^2-1)(N^2-4)(N^2-9) \dots (N^2-(n-1)^2)}.
		\end{equation}

\section{Appendix A: Examples}

	\subsection{A.1:  Class expansion of $m_{\lambda}(\Xi_n)$ for $|\lambda| \leq 4.$}
	
	\begin{description}
	
	\item{$|\lambda|=1$}
\begin{equation*}
m_{(1)}(\Xi_n)= \cc_{(1)}(n).
\end{equation*}

\item{$|\lambda|=2$}
\begin{align*}
m_{(2)}(\Xi_n)=& 
\cc_{(2)}(n)+\frac{1}{2}n(n-1)
\cc_{(0)}(n).  \\
m_{(1^2)}(\Xi_n)=&  
\cc_{(2)}(n) + \cc_{(1^2)}(n). \\
h_2(\Xi_n)=&
2 \cc_{(2)}(n) + \cc_{(1^2)}(n)+\frac{1}{2}n(n-1) \cc_{(0)}(n).
\end{align*}

\item{$|\lambda|=3$}
\begin{align*}
m_{(3)}(\Xi_n)=& 
\cc_{(3)}(n)+(2n-3) \cc_{(1)}(n).\\
m_{(2,1)}(\Xi_n)=& 
3 \cc_{(3)}(n) + \cc_{(2,1)}(n) + \frac{1}{2}(n-2)(n+1)
\cc_{(1)}(n).   \\
m_{(1^3)}(\Xi_n)=& \cc_{(3)}(n) + \cc_{(2,1)}(n) + \cc_{(1^3)}(n). \\
h_3(\Xi_n)=& 
5 \cc_{(3)}(n) + 2 \cc_{(2,1)}(n) + \cc_{(1^3)}(n)
+\frac{1}{2}(n^2+3n-8) \cc_{(1)}(n).
\end{align*}

\item{$|\lambda|=4$}
\begin{align*}
m_{(4)}(\Xi_n)=& \cc_{(4)}(n) +(3n-4) \cc_{(2)}(n)
+4 \cc_{(1^2)}(n) + \frac{1}{6}n(n-1)(4n-5) \cc_{(0)}(n). \\
m_{(3,1)}(\Xi_n) =& 4\cc_{(4)}(n)+ \cc_{(3,1)}(n)+ 2(3n-7) \cc_{(2)}(n)+ 2(2n-3) \cc_{(1^2)}(n) \\
& + \frac{1}{3}n(n-1)(n-2) \cc_{(0)}(n). \\
m_{(2^2)}(\Xi_n)=& 2 \cc_{(4)}(n)+\cc_{(2^2)}(n)+ \frac{1}{2}(n^2-n-4)\cc_{(2)}(n)
+2\cc_{(1^2)}(n)  \\
&+\frac{1}{24}n(n-1)(n-2)(3n-1) \cc_{(0)}(n). \\
m_{(2,1^2)}(\Xi_n) =&
6 \cc_{(4)}(n)+3\cc_{(3,1)}(n)+ 2\cc_{(2^2)}(n)+ \cc_{(2,1^2)}(n)+ \frac{1}{2}(n-3)(n+2) \cc_{(2)}(n) \\
&+\frac{1}{2}(n^2-n-4) \cc_{(1^2)}(n). \\
m_{(1^4)}(\Xi_n)=& \cc_{(4)}(n)+\cc_{(3,1)}(n)+\cc_{(2^2)}(n)+\cc_{(2,1^2)}(n)+\cc_{(1^4)}(n). \\
h_4(\Xi_n)=& 14 \cc_{(4)}(n)+5\cc_{(3,1)}(n)+4\cc_{(2^2)}(n)+2\cc_{(2,1^2)}(n)+\cc_{(1^4)}(n) \\
&+ (n^2+8n-23) \cc_{(2)}(n)+ \frac{1}{2}(n^2+7n-4) \cc_{(1^2)}(n) \\
&\quad + \frac{1}{24}n(n-1)(3n^2+17n-34) \cc_{(0)}(n).
\end{align*}
\end{description}

	\subsection{A.2: Tables of $L_{\mu}^{\lambda}$}
	We now give tables of $L^\lambda_\mu,$ which can be compared with the class
	expansions in Appendix A.1.
	The row labelled ``SUM'' stands for $\sum_{\lambda \vdash k} L^\lambda_\mu$
	for each column associated with $\mu$, which equals $\prod \Cat_{\mu_i}$.	

\begin{center}
\begin{tabular}{|c||c|c|}  \hline
$\lambda \  \backslash \ \mu$ & $2$ & $1^2$  \\ \hline\hline
$2$ & 1 &      \\ \hline
$1^2$ & 1 & 1    \\ \hline \hline
SUM & 2 & 1    \\ \hline
\end{tabular} \quad
\begin{tabular}{|c||c|c|c|}  \hline
$\lambda \  \backslash \ \mu$ & $3$ & $21$ & $1^3$  \\ \hline \hline
$3$ & 1 &  &    \\ \hline
$21$ & 3 & 1 &   \\ \hline
$1^3$ & 1 & 1 & 1  \\ \hline \hline
SUM & 5 & 2 & 1  \\ \hline
\end{tabular} 
\end{center} 
\begin{center}
\begin{tabular}{|c||c|c|c|c|c|} \hline 
$\lambda \  \backslash \ \mu$ & $4$ & $31$ & $2^2$ & $21^2$ & $1^4$ \\ \hline\hline
$4$ & 1 &  &  &  &   \\ \hline
$31$ & 4 & 1 &  &  & \\ \hline
$2^2$ & 2 & 0 & 1 &  &  \\ \hline
$2 1^2$ & 6 & 3 & 2 & 1 &  \\ \hline
$1^4$ & 1 & 1 & 1 & 1 & 1 \\ \hline \hline
SUM & $14$ &  $5$  & $4$ & $2$ & 
$1$ \\ \hline
\end{tabular} \quad
\begin{tabular}{|c||c|c|c|c|c|c|c|} \hline 
$\lambda \  \backslash \ \mu$ & $5$ & $41$ & $32$ & $31^2$ & $2^2 1$ & $2 1^3$ & 
$1^5$ \\ \hline \hline
$5$ & 1 &  &  &  &  &  &    \\ \hline
$41$ & 5 & 1 &  &  &  &  &  \\ \hline
$32$ & 5 & 0  & 1 &  &  &  &  \\ \hline
$3 1^2$ & 10 & 4 & 1 & 1 &  &   &  \\ \hline
$2^21$ & 10 & 2 & 3 & 0 & 1 &  &  \\ \hline
$2  1^3$ & 10 & 6 & 4 & 3 & 2 & 1 &  \\ \hline
$1^5$ & 1 & 1 & 1 & 1 & 1 & 1 & 1 \\ \hline \hline
SUM & $42$ & $14$ & $10$ & $5$ & $4$ & 
$2$ & $1$ \\ \hline
\end{tabular}
\end{center}
\begin{center}
\begin{tabular}{|c||c|c|c|c|c|c|c|c|c|c|c|} \hline 
$\lambda \  \backslash \ \mu$ & $6$ & $51$ & $42$ & $4 1^2$  & $3^2$ 
& $321$ & $3 1^3$ & $2^3$ & $2^2 1^2$ & $2 1^4$ & 
$1^6$ \\ \hline \hline
$6$ & 1 &  &  &  &  &  &  &  &  &  &   \\ \hline
$51$ & 6 & 1 &  &  &  &  &  &  &  &  &  \\ \hline
$42 $& 6 & 0  & 1 &  &  &  &  &  &  &  &  \\ \hline
$4 1^2$ & 15 & 5 & 1 & 1 &  &   &  &  &  &  &  \\ \hline
$3^2$ & 3 & 0 & 0 & 0 & 1 & &  &  &  &  & \\ \hline
$321$ & 30 & 5 & 4 & 0 & 6 & 1 &  &  &  &  & \\ \hline
$3 1^3$ & 20 & 10 & 4 & 4 & 2 & 1 & 1 &  &  &  &  \\ \hline 
$2^3$ & 5 & 0 & 2 &0 & 0 & 0 & 0 & 1 &  &  &  \\ \hline 
$2^2 1^2$  & 30 & 10 & 8 & 2 & 9 & 3 & 0 & 3 & 1 &  &  \\ \hline
$2 1^4$ & 15 & 10 & 7 & 6 & 6 & 4 & 3 & 3 & 2 & 1 &  \\ \hline
$1^6$ & 1 & 1 & 1 & 1 & 1 & 1 & 1 & 1 & 1 & 1& 1 \\ \hline\hline
SUM & $132$ & $42$ & $28$ & $14$ & 25 & 10 & 5 & 8 & 4 & 2 & 1
\\ \hline
\end{tabular}
\end{center}
{\small
\begin{center}
\begin{tabular}{|c||c|c|c|c|c|c|c|c|c|c|c|c|c|c|c|} \hline 
$\lambda \  \backslash \ \mu$ & $7$ & $6 1$ & $5 2$ &  $5 1^2$  & $4 3$ 
& $4 2 1$ & $4 1^3$ & $3^2 1$ & $3 2^2$ & $3 2 1^2$ & $3 1^4$ & $2^3 1$ & $2^2 1^3$ &
$2 1^5$ & $1^7$ \\ \hline \hline
$7$ & $1$ & & & & & & & & & & & & & & \\ \hline
$6 1$ & $7$ & $1$ & & & & & & & & & & & & & \\ \hline
$5 2$ & $7$ & $0$ & $1$ & & & & & & & & & & & & \\ \hline
$5 1^2$ & $21$ & $6$ & $1$ & $1$ & & & & & & & & & & & \\ \hline
$4 3$ & $7$ & $0$ & $0$ & $0$ & $1$ & & & & & & & & & & \\ \hline
$4 2 1$ & $42$ & $6$ & $5$ & $0$ & $3$ & $1$  & & & & & & & & & \\ \hline
$4 1^3$ & $35$ & $15$ & $5$ & $5$ & $1$ & $1$ & $1$ & & & & & & & & \\ \hline
$3^2 1$ & $21$ & $3$ & $0$ & $0$ & $4$ & $0$ & $0$ & $1$ & & & & & & & \\ \hline
$3 2^2$ & $21$ & $0$ & $5$ & $0$ & $2$ & $0$ & $0$ & $0$ & $1$ & & & & & & \\ \hline
$3 2 1^2$ & $105$ & $30$ & $15$ & $5$ & $18$ & $4$ & $0$ & $6$ & $2$ & $1$ & & & & & \\ \hline
$3 1^4$ & $35$ & $20$ & $10$ & $10$ & $5$ & $4$ & $4$ & $2$ & $1$ & $1$ & $1$ & & & & \\ \hline
$2^3 1$ & $35$ & $5$ & $10$ & $0$ & $6$ & $2$ & $0$ & $0$ & $3$ & $0$ & $0$ & $1$ & & & \\ \hline
$2^2 1^3$ & $70$ & $30$ & $20$ & $10$ & $20$ & $8$ & $2$ & $9$ & $7$ & $3$ & $0$ & $3$ & $1$ &  &  \\ \hline
$2 1^5$ & $21$ & $15$ & $11$ & $10$ & $9$ & $7$ & $6$ & $6$ & $5$ & $4$ & $3$ & $3$ & $2$ & $1$ & \\ \hline
$1^7$ & $1$ & $1$ & $1$ & $1$ & $1$ & $1$ & $1$ & $1$ & $1$ & $1$ & $1$ 
& $1$ & $1$ & $1$ & $1$ \\ \hline \hline
SUM & $429$ & $132$ & $84$ & $42$ & $70$ & $28$ & $14$ & $25$ & $20$ & $10$ & $5$ 
& $8$ & $4$ & $2$ & $1$ \\ \hline
\end{tabular}
\end{center}
}

\end{document}